\documentclass[11pt,a4paper]{article}%
\usepackage{ifthen}
\usepackage{amsmath}
\usepackage{amscd}
\usepackage{amssymb}
\usepackage{theorem}
\usepackage{makeidx}%
\usepackage{latexsym}
\usepackage[parfill]{parskip}
\usepackage{graphicx}
\usepackage[english]{babel}

\theoremstyle{plain} 
\newtheorem{thm}{Theorem}[section]
\newtheorem{cor}[thm]{Corollary}
\newtheorem{lem}[thm]{Lemma}
\newtheorem{prop}[thm]{Proposition}
\newtheorem{remark}[thm]{Remark}

\newtheorem{defn}[thm]{Definition}

\theoremstyle{remark}

\newcommand{\thmref}[1]{Theorem~\ref{#1}}
\newcommand{\propref}[1]{Proposition~\ref{#1}}

\newcommand{\lemref}[1]{Lemma~\ref{#1}}
\newcommand{\corref}[1]{Corollary~\ref{#1}}

\newcommand{\defnref}[1]{Definition~\ref{#1}}

\newcommand{\R}{\mathbb{R}}
\newcommand{\C}{\mathbb{C}}

\newcommand{\Log}{\mbox{$\,\mathbb{L}${\rm og}}}
\newcommand{\Trace}{\mbox{$\,\mathbb{T}${\rm race}}}

\newcommand{\Map}{\mbox{{\rm Map}}}

\renewcommand{\=}{\ \approx \ }

\renewcommand{\log}{\textrm{log}\,}

\newcommand{\pdo}{\psi{\rm do}}
\newcommand{\pdos}{\psi{\rm dos}}

\newcommand{\pd}{\partial}

\newcommand{\na}{\nabla}

\newcommand{\noi}{\noindent}

\newcommand{\ob}{{\rm ob}}

\newcommand{\ol}{\overline}
\newcommand{\ul}{\underline}

\newcommand{\To}{\Rightarrow}

\newcommand{\too}{\longrightarrow}

\newcommand{\mto}{\mapsto}

\newcommand{\wh}{\widehat}

\newcommand{\wt}{\widetilde}

\newcommand{\x}{\times}
\newcommand{\ox}{\otimes}
\newcommand{\sox}{\,\mbox{{\small $\ox$}}\,}
\newcommand{\mon}{\mbox{{\tiny $\ox$}}}
\newcommand{\ry}{{\mon y}}
\newcommand{\rw}{{\mon w}}
\newcommand{\lw}{{w \mon }}
\newcommand{\rzz}{{\mon z}}

\newcommand{\ly}{{y \mon}}
\newcommand{\rx}{{\mon x}}
\newcommand{\lx}{{x \mon}}

\newcommand{\aand}{\mbox{and}}

\newcommand{\wwith}{\mbox{with}}

\newcommand{\iin}{\mbox{in}}
\newcommand{\iif}{\mbox{if}}
\newcommand{\iis}{\mbox{is}}
\newcommand{\oon}{\mbox{on}}
\newcommand{\sso}{\mbox{so}}

\renewcommand{\i}{\iota}

\newcommand{\bsh}{\backslash}

\newcommand{\db}{d{\hskip-1pt\bar{}}\hskip1pt}

\newcommand{\Aa}{{\mathcal A}}

\newcommand{\Cc}{{\mathcal C}}
\newcommand{\Dd}{{\mathcal D}}

\newcommand{\Ff}{{\mathcal F}}
\newcommand{\Gg}{{\mathcal G}}

\newcommand{\Kk}{{\mathcal K}}

\newcommand{\Nn}{{\mathcal N}}

\newcommand{\Zz}{{\mathcal Z}}

\newcommand{\as}{{\textsf{a}}}
\newcommand{\Bs}{{\textsf{B}}}
\newcommand{\bs}{{\textsf{b}}}

\newcommand{\Fs}{{\textsf{F}}}

\newcommand{\Gs}{{\textsf{G}}}

\newcommand{\Hs}{{\textsf{H}}}

\newcommand{\Js}{{\textsf{J}}}

\newcommand{\m}{{\textsf{m}}}

\newcommand{\hf}{\hskip 5mm}

\renewcommand{\r}{\right}
\renewcommand{\l}{\left}

\newcommand{\<}{\subset}
\newcommand{\ii}{^{-1}}
\newcommand{\pr}{^{\prime}}

\renewcommand{\d}{\delta}
\newcommand{\Si}{\Sigma}

\renewcommand{\log}{\mbox{{\rm log}\,}}

\newcommand{\z}{\zeta}
\newcommand{\s}{\sigma}

\renewcommand{\a}{\alpha}
\renewcommand{\b}{\beta}

\newcommand{\g}{\gamma}

\newcommand{\la}{\lambda}

\renewcommand{\d}{\delta}

\newcommand{\e}{\varepsilon}

\newcommand{\Z}{\mathbb{Z}}

\renewcommand{\O}{\Omega}

\newcommand{\oo}{\infty}
\newcommand{\mt}{\emptyset}

\newcommand{\End}{\mbox{\rm End}\,}

\newcommand{\Ab}{\textbf{A}}
\newcommand{\Bb}{\textbf{B}}
\newcommand{\Cb}{\textbf{C}}

\newcommand{\Eb}{\textbf{E}}

\newcommand{\Mb}{\textbf{M}}
\newcommand{\Sb}{\textbf{S}}

\newcommand{\Ring}{\textbf{Ring}}
\newcommand{\Abelian}{\textbf{Abelian}}

\newcommand{\Sab}{\mathcal{A}}

\newcommand{\Hom}{\mbox{\rm Hom}}
\newcommand{\Fred}{\mbox{\rm Fred}}

\newcommand{\ind}{\mbox{\rm ind\,}}

\newcommand{\Ind}{\mbox{\rm Ind\,}}

\newcommand{\Ker}{\textmd{\small {\rm Ker}}\,}

\newcommand{\cok}{\mbox{\rm coker}}

\newcommand{\tr}{\mbox{\rm tr\,}}

\newcommand{\Tr}{\mbox{\rm Tr\,}}

\newcommand{\psdo}{{\rm \Psi DO}}

\newcommand\Det{\mbox{\rm Det\,}}

\newcommand{\res}{\textrm{res}\,}

\renewcommand{\t}{\tau}

\newcommand{\sgn}{{\rm sgn}}

\newcommand{\ran}{{\rm ran}}

\newcommand{\cob}{\textbf{Bord}}

\newcommand{\mor}{ \mbox{\textsf{{\rm mor}}}}
\newcommand{\emor}{ \mbox{\textsf{{\rm end}}}}

\newcommand{\alg}{\textbf{Alg}}

\renewcommand{\o}{\bullet}

\renewcommand{\det}{\mbox{{\rm det}\,}}

\numberwithin{equation}{section}

\begin{document} 

\begin{center} 
\Huge{log TQFT}\\[4mm]

{\small Simon Scott}\\[7mm]

\end{center}

The goal here is to put into place an algebraic theory, or rather a categorification, of logarithmic representations and their log-determinant characters.

The motivation for investigating such logarithmic functors is that they  provide a functorial setting for additive invariants arising as generalised Reidemeister torsions on bordism categories. Invariants of this type may be viewed as semi-classical, positioned between genera (classical bordism invariants) and TQFTs (quantum bordism invariants); the former are homomorphisms
 \begin{equation*}
\mu :\Omega_* \to R
\end{equation*}
on the ring $\O_*$ of  bordism classes of closed manifolds, such as the signature of a 4k dimensional manifold, while a TQFT (topological quantum field theory) of dimension $n$ refers to a symmetric monoidal functor 
 \begin{equation*}
Z:\cob_n \to \Bs
\end{equation*}
from  the bordism category $\cob_n$, whose objects are smooth closed (n-1)-dimensional manifolds $M$ and whose morphisms are $n$-dimensional bordisms, to a target symmetric monoidal category $\Bs$.  

The class of semi-classical bordism invariants considered here arise as characters of log-additive simplicial maps
\begin{equation}\label{logfunctor}
\log:\Nn\cob_n \to \Sab
\end{equation}
from the nerve $\Nn\cob_n$ of the bordism category
to a  simplicial set of rings $\Sab$. Such a map \eqref{logfunctor}, called a {\it log-functor}, associates to each bordism $W\in\mor (M_0,M_1)$ between closed manifolds $M_0$ and $M_1$ a logarithm $\log_{\mbox{{\tiny $M_0\sqcup M_1$}}}(W)$ in a ring $ \Fs(M_0\sqcup M_1)\in \Sab$
along with a hierarchy of compatible inclusions 
$$\Fs(M_0\sqcup M_2)$$
$$\downarrow$$ 
\begin{equation}\label{F diag}
 \ \    \ \  \Fs(M_0\sqcup M_1 \sqcup M_2) \ \  \ \    
\end{equation}
$$ \ \  \nearrow  \ \  \hskip 25mm \ \  \nwarrow \ \  $$ 
$$\Fs(M_0\sqcup M_1) \ \   \hskip 35mm \ \  \ \  \ \ \Fs(M_1\sqcup M_2)$$ 

such that when two bordisms $W\in\mor (M_0,M_1), W\pr\in\mor (M_1,M_2)$ are sewn together there is a log-additive identity in  $\Fs(M_0\sqcup M_1 \sqcup M_2)$
\begin{equation} \label{F log} 
\log_{\mbox{{\tiny $M_0\sqcup M_2$}}}(W\cup_{_{\mbox{{\tiny $M_1$}}}} W\pr)  \=  \log_{\mbox{{\tiny $M_0\sqcup M_1$}}}(W)   + \log_{_{\mbox{{\tiny $M_1\sqcup M_2$}}}}( W\pr),
\end{equation}
where $\=$ indicates  equality modulo finite sums of commutators. Neither commutators nor inclusion maps are seen by categorical trace maps $\t_{\mbox{{\tiny $N$}}}: \Fs(N) \to R$ to a commutative ring $R$ and so, irrespective of in which ring it may be convenient to view the logarithm of a bordism $W$, the resulting log-character $\t (\log W):=\t_{\mbox{{\tiny {\bf $M_0\sqcup M_1$}}}} (\log W) \in R$ is invariantly defined.

Characters of log-TQFTs capture a class of semi-local invariants that are of a somewhat more general nature than the local invariants that occur as genera  but which, in view of the log-additive pasting property, must be far simpler and more restricted (possibly more delicate) than the globally determined  invariants of a TQFT.  Such trace-logs include instances of classical Whitehead  and Reidemeister torsions and the topological signature $\s$ and the (relative)  Euler characteristic $\chi$ (note that $\s$ is a genus while $\chi$ is not).  Log-Determinants of this type can arise formally in semi-classical expansions of Feymann path integrals, such as Reidemeister torsion $T_M(a)$ in the stationary phase expansion of Chern-Simons TQFT  $Z_{{\rm cs}}(M) \sim \sum_a c(a)\,\sqrt{T_M(a)}$ over irreducible flat connections \cite{Witten}.  

On the other hand, generalising the classical topological signature $\s$, higher Novikov signatures are additive with respect to gluing \cite{add higher signature} and may be conjectured to be characters of a  log-TQFT  on $\Nn\cob_n$ ranging (following a suggestion by Ryszard Nest) in Hoschchild homology $HH_k(\Aa)$, the case $k=0$ being the subject of this article. 

\section{Logarithmic representations of monoids}

We begin with the notion of a logarithmic representation of a monoid $\Zz$  into a ring  $\Bs = (\Bs, \cdot, +)$. This is defined to be a homomorphism 
\begin{equation}\label{monoid log}
\log:  \Zz \to  \Bs/[\Bs, \Bs],
\end{equation}
where 
\begin{equation}\label{[B,B]}
 [\Bs, \Bs]= \{ \sum_{1\leq j\leq n} [\b_j, \b\pr_j]\ | \ \b_j, \b\pr_j \in \Bs\}
\end{equation}
is the subgroup  of the abelian group $(\Bs, +)$ consisting of finite sums of commutators $[\b_j, \b\pr_j] := \b_j\cdot \b\pr_j  - \b\pr_j\cdot \b_j$ and $ \Bs/[\Bs, \Bs] :=  (\Bs,+)/[\Bs, \Bs]$ is the abelian quotient group. 
For $\mu, \nu\in \Bs$ we may use the notation
\begin{equation}\label{trace equivalent}
\mu\= \nu \ \ \iif \  \mu - \nu \in  [\Bs, \Bs], \hf \sso \  \mu = \nu \ \iin \,  \Bs/[\Bs, \Bs].
\end{equation}  
Thus, one has 
\begin{equation}\label{monoid log 2}
\ \log(ba) = \log \,a + \, \log\, b
\end{equation}
in $\Bs/[\Bs, \Bs]$, where $ba = b\circ a$ is composition in $\Zz$.  A map
$\ell og:  \Zz \to \Bs$ with $$\ell og(ba) = \ell og (b) + \ell og (a) + \sum_j [c_j, c\pr_j]$$
for some $c_j, c\pr_j\in \Bs$, so $\ell og(ba) \= \ell og (b) + \ell og (a)$ in $\Bs$, defines a logarithm, and if the exact sequence $0 \to [\Bs, \Bs] \to \Bs \to \Bs/[\Bs, \Bs] \to 0$ of abelian groups splits then the converse holds. Sums of logs are logs and so form an abelian group
$\Log(\Zz,\Bs):= \Hom(\Zz, \Bs/[\Bs, \Bs]).$
  
A trace on  $\Bs$ with values in a commutative unital ring  $(R, \cdot, +)$ is a homomorphism of abelian groups $\t:  (\Bs, +) \to  (R, +)$ which vanishes on commutators $\t  ([b, b\pr])=0$, so $[\Bs, \Bs]\< \Ker (\t)$. To give $\t$ is equivalent to an abelian group homomorphism $$\wt \t : \Bs/[\Bs, \Bs]  \to R.$$ Sums of traces are traces, forming an abelian group
$\Trace(\Bs, R)$.  A log-character (or logarithmic determinant or trace-log) on $\Zz$  is an evaluation of the canonical pairing
$$\Trace(\Bs, R) \times \Log(\Zz,\Bs) \to \Hom(\Zz, (R, +)), \hf (\t, \, \log) \mto \wt\t\circ\log.$$
Such a character inherits the log-additivity property for $a, b \in \Zz$
\begin{equation}\label{t log add}
\wt\t(\log ba) = \wt\t(\log \,a) + \wt \t(\log\, b) \hf \iin\ R,
\end{equation}
while composition with an exponential map $\e: R \to A^*, \, \e(x+y)=\e(x)\cdot \e(y),$  into the units of  a commutative ring $A$ associates a multiplicative determinant
$a\mto \det a:= e\circ \wt\t \circ \log (a)$.  

For example, let $\Zz = \Fred$ be the monoid of Fredholm operators on a Hilbert space, and  $\Bs = \Ff$ the ideal of finite-rank operators. The map
\begin{equation}\label{log Fred monoid}
\log : \Fred \to  \Ff/[\Ff, \Ff], \hf \log a := \pi([a, p]), 
\end{equation}
where $p\in \Fred$ is any parametrix  for $a$ and $\pi: \Ff \to  \Ff/[\Ff, \Ff]$ the quotient map, is a logarithm, the abstract Fredholm index of $a$, whilst its numeric log-character with respect to the canonical isomorphism $ \Ff/[\Ff, \Ff] \stackrel{\cong}{\to} \C, c\mto \wt\Tr(c),$ defined by the classical trace $\Tr: \Ff \to  \C$ is the usual integer valued Fredholm index $$\wt\Tr(\log a) = \ind a := \dim\ker (a) - \dim\cok (a)$$ and \eqref{t log add} is the classical additivity property of the index $\ind ba = \ind a + \ind b$. Likewise, on continuous families $\Zz = \Map(M, \Fred)$ of Fredholm operators, with  continuous parametrix, parametrized by a manifold $M$, a log-character can be defined by sending  $\as\in \Map(M, \Fred)$  to its index bundle $\log\as := \Ind \as \in K_0(M)$. The top exterior power operation acts as an exponential map on the commutative ring $K_0(M)$ sending $\Ind \as$ to the isomorphism class of the determinant line bundle $\Det \as$ in the group $A\cong  H^2(M, \Z)$ of complex line bundles over $M$, with the log-additivity property 
$ \Ind \bs\as  = \Ind \as + \Ind \bs$ in $K_0(M)$ exponentiating to the canonical multiplicativity property 
$ \Det \bs\as  = \Det \as \ox \Det \bs$ of the determinant line bundle in $A$. (These facts persist to the case of families of Fredholm operators between non-isomorphic bundles, but need to be stated in terms of log-functors on categories.)
 

Similarly, the odd Chern character admits a log-character description as the character of a logarithm  $\log:  \Zz \to  (\Bs, +)/([\Bs, \Bs] + d\Bs)$  to a differential graded ring  $\Bs = (\Bs, d)$, where $[\Bs, \Bs] + d\Bs$ is the abelian subgroup of sums of graded commutators and exact elements $db$ some $b\in B$. The classical Fredholm determinant (arising as the exponentiated character of a  logarithmic representation of the universal cover of the general linear group) and the suspended eta invariant \cite{Me} are particular instances.

On general categories matters are complicated by the fact that the respective logarithms of a pair of composable morphisms will, in general, take values in different rings, and so log-additivity \eqref{monoid log 2} only becomes meaningful within the higher structure \eqref{F diag}, \eqref{F log}. 
 
\section{Logarithmic representations of categories}

All categories will be assumed to be small. Denote the set of morphisms in a category $\Cb$ between objects $x,y \in \ob(\Cb)$ by $\mor_\Cb(x,y)$, or $\mor(x,y)$, and 
$\emor(x) := \mor(x,x)$.   $\Cb$ is  monoidal if it has a bifunctor $\sox: \Cb \x \Cb \to \Cb$  which is associative  with identity object $1=1_\Cb$ up to coherent isomorphism. Any two coherence  isomorphisms between associativity bracketings of an $n$-fold product $x_1\sox x_2 \sox \cdots \sox x_n$ for $x_j \in \ob(\Cb)$ then coincide. To specify for each  $\s\in S_n$ (symmetric group) a permutation isomorphism
\begin{equation}\label{Sn action}
\underbrace{x_1\sox \cdots \sox x_n}_{:= x} \stackrel{s_\s(x)}{\too} \underbrace{x_{\s(1)}\sox \cdots \sox x_{\s(n)}}_{ := x_\s}
\end{equation}
in $\mor_\Cb(x,x_\s)$ a braiding map $b_{w,y}: w\sox y \to y\sox w$  for each $w,y \in \ob(\Cb)$  is assumed with $b_{y,w} = b_{w,y}\ii$, giving $\Cb$ the structure of a symmetric monoidal category: $\sox$ is then commutative up to coherent isomorphism and \eqref{Sn action} is uniquely defined for each associativity bracketing of $x$ and $x_\s$.   A functor $\Fs: \Cb \to \Ab$ out of a monoidal category $\Cb$  will be said to be strict if $\Fs(x_1\sox \cdots \sox x_n)$ is independent of the choice of associativity bracketing of $x_1\sox \cdots \sox x_n$ and if $\Fs$ maps the coherence isomorphisms to identity morphisms in $\Ab$.  
(The assumption that $\Fs$ is strict can be readily dropped provided one keeps track of the isomorphisms  $\Fs((x\sox y)\sox z) \to \Fs(x\sox (y\sox z))$, and so on; essential, for example, for a braided monoidal category).

\begin{lem}\label{mu sigma} For $x = x_1\sox \cdots \sox x_n$ and $\s\in S_n$ one has a canonical isomorphism
\begin{equation}\label{F symmetric 1}
\mu_\s(x)  := \Fs(s_\s(x))  : \Fs(x) \stackrel{ \cong}{\to} \Fs(x_\s),  
\end{equation}
independent of a choice of associativity bracketing of $x$ or $x_\s$, and  satisfying
\begin{equation}\label{F symmetric 2}
\mu_{\s^\prime\circ\s}(x) =  \mu_{\s^\prime}(x_\s)\circ\mu_\s(x).
\end{equation}
\end{lem}

The \emph{product functors} of a monoidal category $\Cb$ are (iterations of) the functors $\Cb\to \Cb$ obtained by holding fixed one of the inputs of the bifunctor $\sox$:  for $y\in \ob(\Cb)$  the right-product functor $\m_\ry  : \Cb \to \Cb$  takes $x\in \ob(\Cb)$ to $x\sox y \in  \ob(\Cb)$ and $\a\in\mor_\Cb(x,z)$ to $\a\sox \iota \in\mor_\Cb(x\sox y ,z\sox y )$, with $\iota$ the identity morphism,  the left-product functor $\m_\lw (x) = w\sox x$ is defined symmetrically. The product functors are not monoidal.

The following construction allows the classical additivity of logarithms to be promoted to a categorical additivity on composed morphisms.

\begin{defn}  Let $\Cb = (\Cb, \sox)$ be a symmetric monoidal category and let $\Cb^* = (\Cb^*, \sox)$ be a groupoid whose objects are those of $\Cb$ and whose morphisms are a specified closed subclass of the isomorphisms of $\Cb$ (containing the coherence  and permutation isomorphisms \eqref{Sn action}). 

A monoidal product representation of the reduced category $\Cb^*$ into an additive category $\Mb$   is a strict functor 
\begin{equation}\label{mpr}
\Fs : \Cb^* \to \Mb
\end{equation}
along with for each $y\in \ob(\Cb)$ a  \textsl{natural transformation of functors}
\begin{equation}\label{eta}
\eta_\ry : \Fs \To \Fs_\ry
\end{equation}
from $\Fs: \Cb^* \to \Mb$ to  $\Fs_\ry := \Fs \circ \m_\ry: \Cb^* \to \Mb$ compatible with $\sox$ and the braiding. (The functor $\Fs$ is not assumed to be monoidal and in general will not be.)\\
\end{defn} 

\begin{lem}\label{pull back MPR}
If $\Sb$ is a symmetric monoidal category, monoidal product  representations pull-back with respect to symmetric monoidal functors $\Js : \Sb^* \to \Cb^*$.  
\end{lem}

$\Fs$ is designed to represent the set of objects of $\Cb$ with its monoidal product, but not necessarily its morphisms. It is, however, sensitive to the permutation isomorphisms of \lemref{mu sigma}, which  intertwine with the covering maps $\eta_\ry$ as follows. 

\begin{lem}\label{mu sigma eta}  Let $y\in \ob(\Cb)$. A monoidal product  representation defines for each $x\in \ob(\Cb)$ a morphism
\begin{equation}\label{eta y x}
\eta_\ry (x)\in \mor_\Mb(\Fs(x), \Fs(x\sox y))
\end{equation}
covering $\m_\ry$ such that for $x$, $x_\s$ as in \eqref{Sn action} 
\begin{equation}\label{eta mu}
\eta_\ry(x_\s) \circ \mu_\s(x) = \mu_{\s\mon 1}(x\sox y) \circ \eta_\ry(x).
\end{equation}
\end{lem}

Proof:\,  A natural transformation $\eta : \Gs \To \Hs$ of functors $\Gs, \Hs : \Ab \to \Bb$ defines for $x\in \ob(\Ab)$ a morphism $\eta(x)\in\mor_\Bb(\Gs(x), \Hs(x))$ with $\eta(z)\circ \Gs(\a) = \Hs(\a)\circ \eta(x)$ for $\a\in\mor_\Ab(x,z)$. Applied to  $\Gs := \Fs$ and $\Hs := \Fs_\ry$, \eqref{eta} gives $\eta_\ry (x) := \eta(x)$ in \eqref{eta y x}. For \eqref{eta mu}, take  $z= x_\s$ and $\a=s_\s(x)\in\mor(x,x_\s)$, so $\eta(z)\circ \Gs(\a) = \eta_\ry (x_\s) \circ \Fs(s_\s(x)) = \eta_\ry (x_\s) \circ  \mu_\s(x)$ while $\Hs(\a)\circ \eta(x) =  \Fs_\ry(s_\s(x))\circ \eta_\ry (x)$ and 
$$\Fs_\ry(s_\s(x)) = \Fs(\m_\ry(s_\s(x))) = \Fs(s_\s(x) \sox \iota_y) = \Fs(s_{\s\mon 1} (x \sox y))  =  \mu_{\s\mon 1}(x \sox y).$$ 
\begin{flushright}
$\Box$
\end{flushright}

In particular, since $\Fs$ is strict there is for each $x\in\ob(\Cb)$ a canonical inclusion
\begin{equation}\label{1 eta}
\eta_x(1) : F(1) \hookrightarrow F(x).
\end{equation}

\emph{Compatibility} of the $\eta_\ry$ with $\sox$ is the requirement $\eta_{\mon(y \mon z)} = \eta_{\mon z}\circ \eta_{\mon y}$, or, more fully,
\begin{equation}\label{eta hom}
\eta_{\mon(y \mon z)}(x) = \eta_{\mon z}(x\sox y)\circ \eta_{\mon y}(x),
\end{equation}
and compatibility with the braiding that
\begin{equation}\label{eta symm}
\eta_{\mon(w \mon z)}(x)  = \mu_{1_x\mon \s_{z,w}}(x\sox z\sox w)\eta_{\mon(z \mon w)}(x) 
\end{equation}
where $1_x\sox \s_{z,w}$ is the permutation which fixes $x$ and swaps $w$ and $z$. 

A monoidal product representation is \emph{injective} if for each $ x\in \ob(\Cb)$ the morphisms $\eta_\ry(x)$   are left-invertible : there is a  
\begin{equation}\label{ejection mor}
\d_\ry(x)\in \mor_\Mb(\Fs(x\sox  y), \Fs( x))  
\end{equation}
with $\d_\ry(x)  \circ \eta_\ry(x)  = i$,  the identity morphism, and satisfying $\d_{\mon z}\circ \d_\ry = \d_{\mon ( z\mon y)}$.

Somewhat more generally, it is useful to combine the above maps to define \emph{insertion morphisms}
for  $x =x_0\sox \cdots \sox x_n$ and $0\leq k \leq n+1$ and $w\in \ob(\Cb)$
\begin{equation}\label{etaw}
\eta^k_w = \eta^k_w (x): \Fs(x_0\sox \cdots \sox x_n) \to \Fs(x_0\sox \cdots \sox x_{k-1}\sox  w\sox x_k  \cdots \sox x_n)
\end{equation}
by 
\begin{equation}\label{etaw defn}
\eta^k_w (x) = \mu_{\s_{k,n+1}}(x\sox w) \circ \eta_\rw(x),
\end{equation}
where $\s_{k,n+1}$ is the permutation $(0,\ldots,n+1) \to  (0,\ldots,k-1, n+1, k, \ldots, n)$.  By {\it fiat}, $\eta_\ry := \eta^{n+1}_y (x)$ and $\eta_\ly := \eta^0_y (x)$. When it is clear what is meant, the superscript $k$ and the domain specifier $(x)$ may be omitted to write $\eta_w$.

For $\ul w= (w_1, \ldots, w_r) \in \ob(\Si(\Cb))$ the iterated insertion morphism
\begin{equation}\label{etaw iterated}
\eta_{\ul w} := \eta_{w_1} \eta_{w_2}\cdots \eta_{w_r} := \eta_{w_1}\circ\cdots\circ \eta_{w_r} : \Fs(x) \to \Fs(x_{\ul w} )
\end{equation}
is unambiguously defined, independently of the ordering of the $\eta_{w_j}$ (in the sense of  \lemref{eta commute}); here, $x= x_0\sox \cdots \sox x_n$ while $x_{\ul w}$ is the monoidal product of the $x_i$ and $w_l$ in a specified order.  If the $\eta_\rw(x)$ are injective then so is \eqref{etaw iterated}: the \emph{ejection morphism} \begin{equation}\label{piwk}
\d^k_w = \d^k_w (x): \Fs(x_w ) \to \Fs(x), \hf \d^k_w (x) = \d_\rw(x)\circ \mu_{\s_{k,n+1}\ii}(x_w),
\end{equation}
for  $x_w = x_0\sox \cdots \sox x_{k-1} \sox w \sox x_{k+1} \sox \cdots \sox x_n$ and $0\leq k \leq n$ and $w\in \ob(\Cb)$ defines a left-inverse for  $\eta^k_w$. The commutation properties are:

\begin{lem}\label{eta commute}
\begin{equation}\label{eta commutation}
\ \eta^l_z \eta^k_w  = 
   \eta^k_w\, \eta^{l-1}_z, \hf k<l,
\end{equation}
\begin{equation}\label{pi commutation}
\ \d^l_w \d^k_z  = 
   \d^{k-1}_z \d^l_w, \hf  k<l, 
\end{equation}
\begin{equation}\label{pi eta commutation}
\ \d^l_w \eta^k_z  = 
\begin{cases}
   \eta^{k-1}_z \d^l_w & {\rm if} \ k<l, \\[1mm]
   \eta^k_z \d^{l-1}_w & {\rm if}  \ k> l,\\[1mm]
    1 & {\rm if}  \ k =l \ {\rm and}  \ w=z.
\end{cases}
\end{equation}
\end{lem}


Proof:\;  Here, 
$\eta^l_z\, \eta^k_w:= \eta^l_z((x\sox w)_{\s_{k, n+1}})\circ \, \eta^k_w(x)$, where $x= x_1\sox \cdots \sox x_n$, and so on. 
The case $\eta^{n+2}_z \, \eta^{n+1}_w  = \eta^{n+1}_w\, \eta^{n+1}_z$ is 
\begin{equation}\label{etaw 3}
\eta_\rzz(x\sox w)\, \eta_\rw(x) = \mu_{1_x\mon \s_{z,w}}(x\sox z\sox w)
\,\eta_\rw(x\sox z) \,\eta_\rzz(x)
\end{equation}
which is a restatement of the compatibility  \eqref{eta hom}, \eqref{eta symm}.  For the general case one has $\eta^l_z\, \eta^k_w  := \mu_{\s_{l, m+2}}((x\sox w)_{\s_{k, m+1}}\sox z)
\,\eta_\rzz((x\sox w)_{\s_{k, m+1}}) \,\mu_{\s_{k, m+1}}(x\sox w) \, \eta_\rw(x)$, by \eqref{etaw defn}.
From \eqref{eta mu},  $\eta_\rzz(x\sox w) \, \mu_{\s_{k, m+1}}(x\sox w)  = \mu_{\s_{k, m+1}\mon 1_z}(x\sox w \sox z) \, \eta_\rzz(x\sox w)$, hence
 \begin{eqnarray}
\eta^l_z\, \eta^k_w  & = & \mu_{\s_{l, m+2}}((x\sox w)_{\s_{k, m+1}}\sox z)\, \mu_{\s_{k, m+1}\mon 1_z}(x\sox w \sox z) 
\,\eta_\rzz(x\sox w) \, \eta_\rw(x) \nonumber \\ 
& \stackrel{\eqref{etaw 3}}{=} & \mu_{\s_{l, m+2}}((x\sox w)_{\s_{k, m+1}}\sox z)\, \mu_{\s_{k, m+1}\mon 1_z}(x\sox w \sox z)  \mu_{1_x\mon \s_{z,w}}(x\sox z\sox w) \nonumber\\
& & \hf\hf\hf\circ\,  \eta_\rw(x\sox z) \, \eta_\rzz (x)\nonumber \\
 & \stackrel{\eqref{F symmetric 2}}{=} & \mu_{\s_{l, m+2} \circ (\s_{k, m+1}\mon 1_z) \circ  (1_x\mon \s_{z,w})}(x\sox z \sox w) 
\,\eta_\rw(x\sox z) \, \eta_\rzz(x).\label{etaw 6}
\end{eqnarray}
The elementary equality 
$ \s_{l, m+2} \circ (\s_{k, m+1}\mon 1_{m+2}) \circ  (1 \mon \s_{m+1,m+2})    = 
\s_{k, m+2} \circ (\s_{l-1, m+1}\mon 1_{m+2})$
 of permutations then yields \eqref{eta commutation}. The other identities follow similarly. 
\begin{flushright}
$\Box$
\end{flushright}

The identities of \lemref{eta commute}  define a (parametrised weakly) simplicial set with  $p$-simplices 
$$\Delta_p = \{ (\xi,  x_0, \ldots, x_{p-1}) \ | \ \xi\in F( x_0\sox \cdots \sox x_{p-1}), x_j\in\ob(\Cb)\} \ \< \ob( \Mb) \x \ob(\Cc^p)$$
with face maps $d_k: \Delta_p \to \Delta_{p-1}, \  (\xi,  x_0, \ldots, x_{p-1}) \mto (\d^k_{x_k}(\xi),  x_0, \ldots, x_{k-1}, x_{k+1}, \ldots, x_{p-1})$, and, for each  $z\in\ob(\Cb)$, degeneracy maps  $$s_k(z): \Delta_p \to \Delta_{p+1}, \  (\xi,  x_0, \ldots, x_{p-1}) \mto  (\eta^k_z(\xi),  x_0, \ldots, x_{k-1}, z, x_k, \ldots, x_{p-1}).$$
It is `weakly' so insofar as the standard simplicial relation `$d_{j+1}s_j(z)=1$'  need not hold.

The morphisms $\d^k_{\ul w}$ are not needed for the development of logarithms, but, when present, they enable more precision in the statement of some logarithm properties. 
 
\textbf{Example:}  \ The fundamental groupoid  $\Pi_{\leq 1}(X)$ of a smooth manifold $X$ is the category whose objects are the points $x$ of $X$ and morphisms are homotopy classes of smooth paths with collared ends, with monoidal product $\sox:= \sqcup$ disjoint union. A $k$-vector bundle $E\to X$ with flat  connection $\nabla$ defines $\Fs_\na: \Pi_{\leq 1}(X) \to \alg_k$ to the category of finite-dimensional $k$-algebras  by
 assigning to $ \ul x = x_1 \sqcup \cdots \sqcup x_n$  the algebra $\Fs_\na(\ul x) =  \End_k (E_{x_1})\oplus \cdots \oplus \End_k(E_{x_n})$ with $E_x$ the fibre of $E$ over $x\in X$ 
 and to  $\ul \g =(\g_1,\dots, \g_n) \in \mor(\ul x,\ul y)$  the canonical isomorphism $\Fs_\na(\ul x)\cong \Fs_\na(\ul y)$ induced by the (invertible) parallel transports $\ss_\na(\g_i)  \in \Hom(E_{x_i}, E_{y_i})$. Here, \eqref{Sn action} is a permutation of the order of the disjoint union $x_1 \sqcup \cdots \sqcup x_n$ and \eqref{F symmetric 1} the corresponding permutation of the matrices $\ss_\na(\g_i)$, while $1$ is the empty set and $\Fs_\na(1) = \{0\}$ the zero algebra and  \eqref{1 eta} the trivial inclusion. The  $\eta_y$ on $\Fs_\na(\ul x)$ are the canonical linear inclusions; in particular, $\eta_{\mon  y}$ is the map $T \mto T\oplus 0 $, while $\d_{\mon y}$ is the corresponding projection map.

\subsection{Tracial monoidal product representations}\label{section tmpr}

On a category of rings $\mathbf{R}$  one has the quotient functor
$
\Pi :  \mathbf{R}\to \mathbf{R}/ [\mathbf{R}, \mathbf{R}] \ \< \ \Abelian,
$
to the category of abelian groups, already used for logarithms on monoids in \S 1, mapping
$ (R, \cdot, +) \in \ob(\mathbf{R})  \ \mto \  (R,+)/ [R,R].$

\begin{defn}  \label{ptmpt}
A monoidal product representation $\Fs$ of a symmetric monoidal category $\Cb$ is said to be pretracial with respect to a background additive category $\Ab$  if the functor $\Fs$ ranges in the category of rings
$$\Fs: \Cb^* \to \textbf{Ring}$$ 
such that for each $x\in \ob(\Cb)$
$$F(x) = \emor_\Ab(\xi_x) $$
for some  unique $\xi_x\in\ob(\Ab)$, and if the insertion morphisms (degeneracy maps)  $\eta_\ry (x)$ of \eqref{eta y x} 
are ring homomorphisms and the $\mu_\s(x)$ of \eqref{F symmetric 1}
with $x = x_1\sox \cdots\sox x_n$ are ring isomorphisms. 
We may indicate this by $\Fs: \Cb^* \to \textbf{Ring}_{\mbox{{\tiny {\bf Add}}}}.$ 

$\Fs$ is said to be injective if the abelian group homomorphisms $\d_\ry(x)$ of \eqref{ejection mor} preserve commutators:  $\d_\ry(x) ([\Fs(x\sox  y),\Fs(x\sox  y)])\< [\Fs( x),\Fs( x)]$. 
\end{defn}
Here, the ring product in $\emor_\Ab(\xi_x)$ is defined by composition of morphisms and the abelian group product by the additive structure on $\Ab$.

\begin{lem}\label{quotient mpr}   
Let $\Fs$ be pretracial and let $\Fs(\Cb^*)$ be the subcategory of $\Ring_{\mbox{{\tiny {\bf Add}}}}$ with objects $\Fs(x)$ for $x\in \ob(\Cb)$. By composing with the quotient functor,
 $\Fs$ pushes-down to an induced
monoidal product representation
\begin{equation}\label{q mpr}
\Fs_\Pi : \Cb^* \to \Fs(\Cb^*)/[\Fs(\Cb^*), \Fs(\Cb^*)], \hf x\mto \Fs( x)/ [\Fs( x),\Fs( x)].
\end{equation}
\end{lem}

Proof:\;   
Since $\Fs$ is pretracial  
$\eta_{\ul w} : \Fs(x) \to \Fs(x_{\ul w} )$  is a ring homomorphism,  taking commutators to commutators. As such, it pushes-down to a homomorphism of abelian groups
\begin{equation}\label{eta tilde}
\wt \eta_{\ul w} : \Fs(x)/[ \Fs(x),  \Fs(x)] \to \Fs(x_{\ul w})/[ \Fs(x_{\ul w}),  \Fs(x_{\ul w})], \hf 
\wt \eta_{\ul w} ([\xi]) := \pi_x \circ  \eta_{\ul w} (\xi), 
\end{equation}
with
$\pi_x: \Fs(x) \to \Fs(x)/[ \Fs(x),  \Fs(x)] $
the quotient map, defining the insertion maps of a monoidal product representation.  Since \eqref{eta commutation} persists to the quotient, $$(\Fs(\Cb^*)/[\Fs(\Cb^*), \Fs(\Cb^*)] , \wt\eta^j_z)$$ inherits the structure of a presimplicial set, while if $\Fs$ is injective then it inherits the structure of a simplicial set from $\Fs(\Cb^*)$.
$\hfill \Box$

A  monoidal category $\Eb$ has a trace $\t$ if there exist  objects $x\in\ob(\Eb)$ with a non-empty closed subclass $\emor^{\,\t}_\Eb(x)$ of endomorphisms and a  map
$$\t_x: \emor^{\,\t}_\Eb(x) \to \emor_\Eb(1)$$
with the trace property that for $\a\in \mor_\Eb(x,y)$ and $\b\in \mor_\Eb(y,x)$ with $\b\circ\a\in \emor^{\,\t}_\Eb(x) $ and  $\a\circ\b\in \emor^{\,\t}_\Eb(y) $  one has $\t_x(b\circ\a) = \t_y(a\circ\b)  \in  \emor_\Eb(1).$
An element $\d\in\emor^{\,\t}_\Eb(x)$ is called  $\t$-trace class and $\t$ a categorical trace. For example, in $\cob_n$ all bordisms   are trace class for the trace sending $W\in \emor(M)$ to the closed manifold  formed by gluing the two boundary portions $\ol M$ and $M$ of $W$ via  the diffeomorphism $ \pd W \stackrel{\cong}{\to} \ol M \sqcup M$, see  \cite{PoSh}, \cite{StTe}. On the other hand, for the classical trace $\Tr$ on the category  of  Hilbert spaces only preferred sub ideals of bounded operators are trace class. Nevertheless, the $\t$ superscript in $\emor^{\,\t}_\Eb(x)$ will be omitted with the understanding that, where necessary, statements are meant for trace class morphisms. 



\begin{defn}  \label{TMPR}
A pre-tracial monoidal product representation  $\Fs: \Cb^* \to \Ring_{\mbox{{\tiny {\bf Add}}}}$ is said to be a tracial monoidal product representation of $\Cb$ if $\Ab$ has an $\Fs$-compatible trace $\t$. $\Fs$-compatible means that $\t$ assigns to each $x\in\ob(\Cb)$ a trace $\t_x : \Fs(x) = \emor_\Ab(\xi_x) \to \emor_\Ab(1_\Ab)$   satisfying the compatibility requirement that for all $x,y \in\ob(\Cb)$ 
\begin{equation}\label{trace eta}
\t_{x\mon y} \circ \eta_\ry(x) = \t_x \hf \aand\hf \t_{x_\s} \circ \mu_\s(x)  = \t_x.
\end{equation}
\end{defn}

Characters in a tracial monoidal product representation can be computed `anywhere': 

\begin{lem}\label{two traces}   
For a tracial monoidal product representation one has
\begin{equation}\label{two traces 2}
\t_{x} = \t_{x_w} \circ \eta_{\ul w}.
\end{equation}
\end{lem}

Proof:\;  Replacing  $\t_{x\mon  z} $ by $\t_{x\mon w\mon  z} \circ \eta_w$  defines another trace on on $\Fs(x\sox z)$, but
$$\t_{x\mon w\mon  z} \circ \eta_w \stackrel{\eqref{eta symm}}{=}  \t_{x\mon w\mon  z} \circ \mu_\s(x\sox z\sox w)\circ \eta_\rw(x\sox z) \stackrel{\eqref{trace eta}}{=} \t_{x\mon z\mon  w} \circ  \eta_\rw(x\sox z) \stackrel{\eqref{trace eta}}{=} \t_{x\mon  z}.$$
Then \eqref{two traces 2} follows by iteration.
$\hfill\Box$

Each of the above structures pushes-down to the quotient monoidal product representation $\Fs_\Pi$ (noted in \eqref{eta tilde} for the insertion maps) while for the trace $\t$ one has for each object $x\in\ob(\Cb)$ a commutative diagram 
\begin{equation*}
\begin{array}{ccl}
  \Fs(x) \ &  \stackrel{\t_x}{\rightarrow} &   \emor_\Eb(1)\\[3mm]
  \downarrow\pi_x  &  \stackrel{\tilde\t_x}{\hf\nearrow} & \\[3mm]
\frac{\Fs(x)}{[\Fs(x), \Fs(x)]}   & &   \end{array}.
\end{equation*}
From this view point, $\pi_x$ is a `universal trace' on $\Fs(x)$ insofar as any trace  factors uniquely through it: one has
$
\t_{x} = \wt\t_{x} \circ \pi_x $ and $ \wt \t_{x} = \wt\t_{x_{\ul w}} \circ \wt \eta_{\ul w},$
with the second identity consequent on \eqref{two traces 2}. Matters may be summarised as the commutativity of the diagram

\begin{equation}\label{trace diagram} 
\begin{array}{rcccl}
  \Fs(x) \ &  & \stackrel{ \eta_{\ul w}}{\longrightarrow} &  & \  \Fs(x_{\ul w})   \\
   & \stackrel{\t_x}{\searrow} & \hf \hf\hf\hf &  \stackrel{\t_{x_{\ul w}}}{\hf\swarrow}  \hskip 5mm &    \\
    &  & &  &   \\
  \downarrow \pi_x  &  & \C &  & \hf\downarrow\pi_{x_{\ul w}} \\
    &  & &  &  \\
  & \stackrel{\tilde\t_x}{\hf\nearrow} \hskip 5mm&  & \stackrel{\tilde\t_{x_{\ul w}}}{\nwarrow} &    \\
 \frac{\Fs(x)}{[\Fs(x), \Fs(x)]}  &  & \stackrel{\tilde\eta_{\ul w}}{\longrightarrow} &  & \frac{\Fs(x_{\ul w})}{[\Fs(x_{\ul w}), \Fs(x_{\ul w})]} .\\[2mm] 
\end{array}
\end{equation}

In particular, (repeating \eqref{eta tilde})
$
 \pi_{x_{\ul w}} \circ  \eta_{\ul w}  = \wt\eta_{x_{\ul w}} \circ  \pi_x.$

\subsection{Logarithmic functors} 

The nerve $\Nn\Cb$ of a category $\Cb$ is the simplicial set whose $p$-simplices are diagrams 
\begin{equation}\label{psimplex}
x_0\stackrel{\a_0}{\to} x_1 \stackrel{\a_1}{\to} x_2 \to \cdots \to x_{p-1} \stackrel{\a_{p-1}}{\to} x_p \hf\in \Nn_p\Cb
\end{equation}
of  morphisms $\a_j \in\mor(x_j, x_{j+1})$.  The $j^{{\tiny{\rm th}}}$ face map 
$d_j : \Nn_p\Cb \to \Nn_{p-1}\Cb$  of the simplex deletes $x_j$, replacing when $0<j<p$
\begin{equation*}
\cdots \to x_{j-1}\stackrel{\a_{j-1}}{\to} x_j \stackrel{\a_j}{\to} x_{j+1} \to \cdots \hf {\rm by} \hf \cdots \to x_{j-1}\stackrel{\a_j\circ \a_{j-1}} {\too}x_{j+1} \to \cdots
\end{equation*}
and  the $j^{{\tiny{\rm th}}}$ degeneracy map $s_j : \Nn_p\Cb \to \Nn_{p+1}\Cb$ replaces 
\begin{equation}\label{sj}
\cdots \to x_j  \stackrel{\a_j}{\to} x_{j+1} \to \cdots \hf {\rm by}\hf \cdots \to x_j \stackrel{\i}{\to} x_j   \stackrel{\a_j}{\to} x_{j+1} \to \cdots .
\end{equation}
$\Nn\Cb$ carries more data than $\Cb$ ---  the objects and morphisms of $\Cb$  are respectively identified with $\Nn_0\Cb$ and  $\Nn_1\Cb$, while there is no right inverse to the composition face map $d_1: \mor_{x_1}(x_0, x_2) \to \mor(x_0, x_2)$. The classifying space  $B\Cb$ of $\Cb$ is  the geometric realisation of $\Nn\Cb$.  

Logarithms on a category $\Cb$ have to be differentiated between according to the substrata of marked morphisms in $\Nn_p\Cb$ on which they act. To this end, one has the  stratum of $\ul z = (x_1,  \ldots, x_{p-1})$-marked $p$-simplices \eqref{psimplex} between $x, y \in\ob(\Cb)$   
$$\mor_{\ul z}(x,y) = \{ x\stackrel{\a_0}{\to} x_1 \stackrel{\a_1}{\to} x_2 \to \cdots \to x_{p-1} \stackrel{\a_{p-1}}{\to} y\} \ \< \  \Nn_p\Cb$$
$$ :\cong \mor_\Cb(x, x_1) \x \mor_\Cb(x_1, x_2) \x \cdots\x \mor_\Cb(x_{p-1}, y).$$
If $\mor(x_j, x_{j+1})=\emptyset$ some $j$ then $\mor_{\ul z}(x,y) :=\emptyset$, while $\mor_\emptyset(x,y) := \mor(x,y)$.    One has  the composition 
\begin{equation*}
 \mor_{\ul z}(x,w)\x  \mor_{\ul z^\prime}(w,y) \stackrel{\circ}{\to} \mor_{\ul z \o w \o \ul z^\prime} (x,y), 
\end{equation*}
relative to concatenation $\o$, so $(x,z)\o y = (x,z,y)$ and so on, as a partially defined composition  $$\Nn_p\Cb\x\Nn_q\Cb\to \Nn_{p+q-1}\Cb$$ on compatible strata, while the face and degeneracy maps respectively restrict to simplicial maps
\begin{equation*}
d_j : \mor_{\ul z}(x,y) \to \mor_{\d_j(\ul z)}(x,y),  \hf s_j : \mor_{\ul z}(x,y) \to \mor_{\s_j(\ul z)}(x,y)
\end{equation*}
with $\d_j : \Cb^p \to \Cb^{p-1}$ and $\s_j : \Cb^p \to \Cb^{p+1}$ defined in the evident way. 

Recall that a simplicial map $f: X \to X\pr$ between simplicial sets $(X, d_j, s_j), (X\pr, d\pr_j, s\pr_j)$ is given by  maps $f_p : \Delta_p  \to \Delta\pr_p$ between $p$-simplices which commute with the face and degeneracy maps, so that
$f_{p-1} d_j = d\pr_j f_p$ and $f_p s_j = s\pr_j f_{p-1}.$
Both these are implied by (but do not imply) 
\begin{equation}\label{sj}
s\pr_j f_{p-1}d_j =  f_p.
\end{equation}
\eqref{sj} is advantageous, here, insofar as it does not involve the boundary operators $d\pr_j$ on $X\pr$. In the case where the range is only a presimplicial set $(X\pr, s\pr_j)$, so that $s\pr_ l s\pr_k = s\pr_k s\pr_{l-1}$ for $k<l$,  a map $f:(X, d_j, s_j) \to (X\pr, s\pr_j)$ may be said to be {\em presimplicial} if  \eqref{sj}  holds. (This applies equally when the domain is also only presimplicial $(X, d_j)$.) 

\begin{defn}\label{log representation}
Let $\Cb = (\Cb, \ox)$ be a symmetric monoidal category and let $$\Fs: \Cb^* \to  {\rm {\bf Ring}}_{\mbox{{\tiny {\bf Add}}}}$$ be a (strict) pretracial monoidal product  representation.  Then a log-functor (or logarithmic-functor) on $\Cb$ taking values in $\Fs$ is a presimplicial log-additive map 
\begin{equation}\label{log functor}
\log: (\Nn\Cb, d_j, s_j)  \to  (\Fs(\Cb^*)/ [\Fs(\Cb^*), \Fs(\Cb^*)], \wt\eta^{\,j}).
\end{equation}
Such a structure is said to define a logarithmic representation of  $\Cb$.
\end{defn}

Unwrapping the definition, a log-functor comprises the following:

1. \  A (strict) pre-tracial monoidal product  representation (on the set $\Nn_0\Cb$ of $0$-simplices): 
$\Fs: \Cb^* \to  \Ring_{\mbox{{\tiny {\bf Add}}}}$, and hence a quotient monoidal product representation $$\Cb^* \to  \Fs(\Cb^*)/[\Fs(\Cb^*), \Fs(\Cb^*)], \hf z\in \ob(\Cb) \mto \Fs(z)/[\Fs(z), \Fs(z)],$$
with insertion maps
$$\wt\eta_{\ul w}:\Fs(z)/[\Fs(z), \Fs(z)] \to \Fs(z_{\ul w})/[\Fs(z_{\ul w}), \Fs(z_{\ul w})].$$

2. \  A simplicial system of (strict) logarithm maps  (on the set $\Nn_1\Cb $ of $1$-simplices) assigning to $x, y \in \ob(\Cb)$, with $x,y$ not both the monoidal identity $1\in \ob(\Cb)$, a map
\begin{equation}\label{logxy}
\log_{x\mon y} : \mor(x,y) \to \Fs(x \sox y)/[\Fs(x \sox y), \Fs(x \sox y)],
\end{equation}
$$ \a\mto \log_{x\mon y}\a = \log(x \stackrel{\a}{\to} y)$$
and, more generally, (on the set $\Nn_p\Cb$ of $p$-simplices) to each marking $\ul z = (z_1,  \ldots, z_{p-1})$  a map

\begin{equation}\label{logxzy}
\log_{x\mon\ul z\mon y} : \mor_{\ul z}(x,y) \to\Fs(x\sox \ul z\sox y)/[\Fs(x\sox \ul z\sox y), \Fs(x\sox \ul z\sox y)] 
\end{equation}
where $ x\sox \ul z\sox y  := x\sox z_1\sox \cdots \sox z_{p-1}\sox y \neq 1$, 
$$\ul \a\mto \log_{x\mon\ul z\mon y}\,\ul \a := \log_{x\mon\ul z\mon y}(x\stackrel{\a_0}{\to} z_1 \stackrel{\a_1}{\to} z_2 \to \cdots \to z_{p-1} \stackrel{\a_{p-1}}{\to} y ),$$
such that for $x \stackrel{\a}{\to} z \stackrel{\b}{\to} y \in \mor_z(x,y)$ associated to $\a\in \mor(x,z)$ and $\b\in \mor(z,y)$ one has in
\begin{equation}\label{quotient F}
\Fs(x\sox z\sox y)/[\Fs(x\sox z\sox y), \Fs(x\sox z\sox y)]
\end{equation}
the ($p=2$)  log-additive property
\begin{equation}\label{log equality} 
\log_{x\mon z\mon y} (x \stackrel{\a}{\to} z \stackrel{\b}{\to} y)  := \wt\eta_\ry(\log_{x\mon z} \,\a) + \wt\eta_\lx(\log_{z\mon y}\, \b),
\end{equation}
or, equivalently, 
\begin{equation}\label{log equality 2} 
\wt\eta_z(\log_{x\mon y} (x \stackrel{\b\circ\a}{\to} y) )  = \wt\eta_\ry(\log_{x\mon z} \,\a) + \wt\eta_\lx(\log_{z\mon y}\, \b).
\end{equation}

{\it Notation}: For brevity, in the left-hand side of \eqref{log equality} and \eqref{log equality 2} we write 
$$\log_{x\mon z\mon y} \,\b \a : = \log_{x\mon z\mon y} (x \stackrel{\a}{\to} z \stackrel{\b}{\to} y), \hf
\log_{x\mon y} \,\b \a : = \log_{x\mon y} (x \stackrel{\b\circ\a}{\to} y).$$

In practise, \eqref{log equality}  is generally obtained consequent on an equivalence
$$\log_{x\mon z\mon y} \,\b\a  =  \eta_\ry(\log_{x\mon z} \,\a) + \eta_\lx(\log_{z\mon y} \,\b) + \sum_{1\leq j\leq m} [\nu_j, \nu^\prime_j]$$ 
some  $\nu_j, \nu^\prime_j\in \Fs(x\sox z\sox y)$ and, 
likewise for \eqref{log equality 2}. 
In this case, the presimpliciality of the log maps \eqref{logxy}, \eqref{logxzy} is for $p=2$
\begin{equation}\label{log compatible}
\log_{x\mon z\mon y}  (x \stackrel{\a}{\to} z \stackrel{\b}{\to} y)  -   \eta_z \,\log_{x\mon y} (x \stackrel{\b\circ\a}{\to} y) \ \in \ [\Fs(x\sox  z\sox y), \Fs(x\sox z\sox y)] 
\end{equation}
\begin{equation}\label{log compatible 2}
\log_{x\mon x\mon y}  (x \stackrel{\a}{\to} x \stackrel{\b}{\to} y)  -   \eta_\lx \,\log_{x\mon y} (x \stackrel{\b\circ\a}{\to} y) \ \in \ [\Fs(x\sox  x\sox y), \Fs(x\sox x\sox y)] 
\end{equation}
\begin{equation}\label{log compatible 3}
\log_{x\mon x\mon y}  (x \stackrel{\a}{\to} y \stackrel{\b}{\to} y)  -   \eta_\ry \,\log_{x\mon y} (x \stackrel{\b\circ\a}{\to} y) \ \in \ [\Fs(x\sox  y\sox y), \Fs(x\sox y\sox y)] 
\end{equation}
and, more generally, with $\ul z = (x_1,  \ldots, x_{p-1})$,  $\nu \in\mor_{\ul z}(x,y)$, $j\in\{1, \dots, p-1\},$ that
\begin{equation}\label{log compatible z}
\log_{\ul z} \, \nu  -   \eta_{x_j}(\log_{\d_j(\ul z)} \,d_j(\nu)) \ \in \ [\Fs(x\sox\ul z\sox y), \Fs(x\sox\ul z\sox y)]
\end{equation}
plus the corresponding two end-point special cases ($x_0=x, x_p=y$) generalising \eqref{log compatible 2} and \eqref{log compatible 3}. These are the identities \eqref{sj} for the presimplicial structures at hand. 

\begin{remark}

{\rm [1] \ A log-functor is not in general a functor of categories, but is a functor of $\oo$-categories.

[2] \  Taking the geometric realization of (both sides of) \eqref{log functor} gives a `logarithm' representation
$|\log| : B\Cb \to  |(\Fs(\Cb^*)/ [\Fs(\Cb^*), \Fs(\Cb^*)]|$ of the (pre-) classifying space $B\Cb$ of the category $\Cb$.} 
\end{remark}

The intertwining of the logarithm and the simplicial structures is clear when written as:
\begin{lem}
The log-additivity property  \eqref{log equality 2}  can be written
\begin{equation*}\label{log equality simplicial}
\wt\eta_1\,\log_{\d_1(\ul x)} \l(d_1(x \stackrel{\a}{\to} z \stackrel{\b}{\to} y)\r)  =  \wt\eta_0\, \log_{\d_0(\ul x)} \l(d_0(x \stackrel{\a}{\to} z \stackrel{\b}{\to} y)\r) + \wt\eta_2\,\log_{\d_2(\ul x)} \l(d_2(x \stackrel{\a}{\to} z \stackrel{\b}{\to} y)\r).
\end{equation*}
where $\ul x = x\sox y\sox z$, $\eta_0 := \eta_\lx, \eta_1 := \eta_z, \eta_2 := \eta_\ry$, $ x \stackrel{\a}{\to} z \stackrel{\b}{\to} y \in \mor_z(x,y)\in \Nn_2\Cb$. 
\end{lem}
Here, the end-point face maps $d_0, d_p: \Nn_p\Cb \to \Nn_{p-1}\Cb$ are defined by deleting the $0^{{\tiny {\rm th}}}$ or  $p^{{\tiny {\rm th}}}$ morphism from a simplex; and  the reason that  \eqref{log compatible 2}, \eqref{log compatible 3} are stated separately.

We note that a log-functor is effectively determined by its action on 1-simplices:

\begin{lem}\label{1-simplices rule ok}  A simplicial system  of logarithm maps $\log_{x\mon \ul z\mon y}$ is determined up to terms in $[\Fs, \Fs]$ \textsl{by}  the log maps $\log_{x\mon y}$ on $\mor(x,y)$ for each $x,y\in\ob(\Cb)$. To define a compatible system of logarithm maps $\log_{x\mon \ul z\mon y}$  it is enough to define the $\log_{x\mon y}$ \textsl{on} $\mor(x,y)$ satisfying \eqref{log equality 2}.
\end{lem}
Proof:\;   Compatibility \eqref{log compatible} gives
$
\log_{x\mon \ul z\mon y} \,\d = \wt\eta_{\ul z}(\log_{x\mon y} \,\d)$ in $\Fs(x\sox  \ul  z\sox y)/[\Fs(x\sox  \ul  z\sox y), \Fs(x\sox  \ul  z\sox y)]$
which is the first statement of the lemma.
Given $\log_{x\mon y}$,  the second statement is that
$
\log_{x\mon \ul z\mon y} \,\d := \wt\eta_{\ul z}(\log_{x\mon y} \,\d),
$
defines  by default  a compatible system of logs \eqref{logxzy}. 
\begin{flushright}
$\Box$
\end{flushright} 

Two $p$ simplices which collapse to the same $(p-r)$ simplex have the same logarithm, and, likewise, inflating simplices does not change logarithms:
\begin{lem}
If  $d_1(x\stackrel{\a}{\to} z \stackrel{\b}{\to} y)  = d_1(x\stackrel{\a\pr}{\to} z \stackrel{\b\pr}{\to} y)$ (that is, $\b\a  = \b\pr\a\pr$) in $\mor(x,y)$ then 
\begin{equation}\label{log well-defined}
\log_{x\mon z\mon y} \,\b\a =  \log_{x\mon z\mon y} \,\b\pr\a\pr
\end{equation}
in $\Fs(x\sox z\sox y)/[\Fs(x\sox z\sox y), \Fs(x\sox z\sox y)] $. More generally, if for $\ul z = (x_1,  \ldots, x_{p-1})$  and $\nu, \nu\pr \in\mor_{\ul z}(x,y)$ and $j=1, \dots, p-1$ one has $d_j(\nu) = d_j(\nu\pr)$, then 
\begin{equation}\label{log collapsing simplices}
\log_{\ul z} \, \nu  = \log_{\ul z} \, \nu\pr
\end{equation}
in $\Fs(x\sox\ul z\sox y)/[\Fs(x\sox\ul z\sox y), \Fs(x\sox\ul z\sox y)]$. Iteratively, if $d_k(d_j(\nu)) = d_k(d_j(\nu\pr))$ then \eqref{log collapsing simplices} continues to hold  since 
\begin{equation}\label{log collapsing two simplices}
\log_{\ul z} \, \nu  =   \wt\eta_{x_j} \wt\eta_{x_k}\,\log_{\d_k(\d_j(\ul z))} \, d_k(d_j(\nu)).
\end{equation}
For $j<k$ 
\begin{equation}\label{commute collapsed simplices}
\eta_{x_j} \eta_{x_k}\,\log_{\d_k(\d_j(\ul z))} \, d_k(d_j(\nu)) = \wt\eta_{x_{k+1}} \wt\eta_{x_j}\,\log_{\d_j(\d_{k-1}(\ul z))} \, d_j(d_{k-1}(\nu)).
\end{equation}
Dually, for the degeneracy maps \eqref{sj} one has 
\begin{equation}\label{log inflating two simplices}
\log_{\s_j(\ul z)} \, s_j(\nu)  =   \wt\eta^j_{x_j} \log_{\ul z} \,\nu 
\end{equation}
\begin{equation}\label{log inflating two simplices 2}
\log_{\s_k(\s_j(\ul z))} \, s_k(s_j(\nu))  =   \wt\eta^k_{x_k} \eta^j_{x_j} \log_{\ul z} \, \nu
\end{equation}
and a corresponding commutation formula to \eqref{commute collapsed simplices}.  For each of the above, the two end-point special cases corresponding to  \eqref{log compatible 2} and \eqref{log compatible 3} also hold.
\end{lem}

Proof:\;   By \eqref{log compatible}
$$ \log_{x\mon z\mon y} (x\stackrel{\a}{\to} z \stackrel{\b}{\to} y) = \wt\eta_z \log_{x\mon y} (x\stackrel{\b\a}{\to} y) = \wt\eta_z \log_{x\mon y} (x\stackrel{\b\pr\a\pr}{\to} y) = \log_{x\mon z\mon y} (x\stackrel{\a\pr}{\to} z \stackrel{\b\pr}{\to} y),$$
and in general $\log_{\ul z} \, \nu  =  \wt\eta_{x_j}(\log_{\d_j(\ul z)} \,d_j(\nu))  =  \wt\eta_{x_j}(\log_{\d_j(\ul z)} \,d_j(\nu\pr)) = \log_{\ul z} \, \nu$ by \eqref{log compatible z}. The general version follows by iterating these equalities given that  \eqref{log collapsing two simplices} holds, and that holds because the $\eta_{x_l}$ are ring homomorphisms. \eqref{commute collapsed simplices} and its $s_j$ counterpart are immediate from \eqref{eta commute} and the simplicial identities $d^j d^k = d^k d^{j-1}$ and $s^j s^k = s^k s^{j+1}$ for $k<j$.  The inflation formulae \eqref{log inflating two simplices}, \eqref{log inflating two simplices 2} follow from \eqref{log compatible z} (resp. \eqref{commute collapsed simplices}) by  replacing $\nu$ by $s_j(\nu)$ (resp. $s_k(s_j(\nu))$). The two end-point special cases of \eqref{log collapsing simplices}
hold from \eqref{log compatible 2} and \eqref{log compatible 3} by the same argument as the case $1\leq j\leq p-1$, while for \eqref{log inflating two simplices} this is shown in \propref{log properties} (2.). 
\begin{flushright}
$\Box$
\end{flushright}

Log-functors transform naturally: if $\Js : \Sb \to \Cb$ is a symmetric monoidal functor, then, since $\Cb \to \Nn\Cb$ is functorial, a logarithmic representation of $\Cb$ pulls-back to one of $\Sb$. Further basic properties of log-functors are listed in the following lemma:

\begin{prop}\label{log properties} 1. Let $p\in \mor_\Cb(x,x)$ be a projection morphism: $p\circ p = p$. Then in $\Fs(x\sox x\sox x)$
\begin{equation}\label{log eq gen}
\eta_\lx (\log_{x\mon x} \,p) \=0. 
\end{equation}
In particular, $\eta_\lx (\log_{x\mon x} \,\iota) \=0$,  where $\iota$ is the identity morphism. 
If $\Fs$ is injective, in the sense of \defnref{TMPR}, then in $\Fs(x\sox x)$
\begin{equation}\label{log eq gen 2}
\log_{x\mon x} \,p \=0. 
\end{equation}

2. For $\a \in \mor(x,y)$ and identity morphisms $\iota_x \in \mor(x,x)$, $\iota_y \in \mor(y,y)$
\begin{equation}\label{log ia}
\log_{x\mon y \mon y} \,(\iota_y \circ \a) \= \eta_\ry (\log_{x\mon y} \a)  \hf\hf \mbox{\textsl{in}} \ \,\Fs(x\sox y \sox y),
\end{equation}
\begin{equation}\label{log ai}
\log_{x\mon x \mon y} \,(\a \circ \iota_x)  \= \eta_\lx (\log_{x\mon y} \a)  \hf\hf \mbox{\textsl{in}} \ \,\Fs(x\sox x \sox y).
\end{equation}
Notation: $\log_{x\mon y \mon y} \,(\iota_y \circ \a): = \log_{x\mon y \mon y} \,(x\stackrel{\a}{\to} y \stackrel{\iota_y}{\to} y)$ .

3.  \textsl{For} $\a, \b \in \mor(x,x)$  one has \textsl{in} $\Fs(x\sox x \sox x)$
\begin{equation}\label{log ba}
\eta_\rx\log_{x\mon x} \,\b\a \= \eta_\rx \log_{x\mon x} \,\a +  \eta_\rx\log_{x\mon x}\, \b. 
\end{equation}

4. For $\a \in \mor(x,x)$ and an isomorphism 
$q\in \mor(w,x)$ one has in $\Fs(w\sox x\sox x\sox w)$
\begin{equation}\label{log qii a q}
\log_{w\mon x\mon x \mon w} \,(q\ii \a q) \=  \eta_\lw\eta_\rw(\log_{x\mon x} \,\a). 
\end{equation}
In the case $x=w$, considering $q\ii \a q\in\mor(x,x)$, if $\Fs$ is injective then
\begin{equation}\label{log qii a q x=w}
\log_{x\mon x } \,(q\ii \a q) \=  \log_{x\mon x} \,\a
\end{equation}
in $\Fs(x\sox x)$. In either case, for a log-determinant structure one has in $\mor_\Ab(1,1)$
\begin{equation}\label{trace log qii a q}
\t(\log \,q\ii \a q) =  \t(\log \,\a)
\end{equation}
for any choice of representatives $\log_{x\mon \ul w\mon x} \, q\ii \a q$ and $\log_{x\mon \ul w\mon x} \a$ of the logarithms.
 
5. Let
$\ul w, \ul w^\prime \in \ob(\Si(\Cb))$ and let
 $\a \in \mor_{\ul w}(x,z)\<\Nn_p\Cb, \b\in \mor_{\ul w^\prime }(z,y)\<\Nn_q\Cb$. Then for a logarithmic representation one has in $\Fs(x\sox  \ul w \sox z \sox \ul w^\prime  \sox y)$ 
\begin{equation}\label{log equality gen}
\log_{x\mon \ul w \mon z\mon \ul w^\prime \mon y} (\b\a)  \=  \eta_{\ul w^\prime \o y}(\log_{x\mon \ul w \mon z} \,\a) + \eta_{x\o \ul w}(\log_{z\mon \ul w^\prime \mon y} \,\b).
\end{equation}

6. Let $\ul w=(w_1,\dots, w_m)\in\ob(\Si(\Cb))$ and let  $\a  = \a_{m+1} \a_m\cdots\a_1\in \mor_{\ul w}(x,y)$ with $\a_j: w_{j-1} \to w_j$ and $w_0:=x$, $w_{m+1}:=y$. Then 
\begin{eqnarray*}\label{log a1am}
\eta_{\ul w}\,\log_{x \mon  y} (\a_{m+1} \a_m\cdots \a_1)  = \log_{x\mon \ul w \mon  y} (\a_{m+1} \a_m\cdots \a_1) =  \sum_{j=1}^{m+1} \eta_{j-1,j}\l(\log_{w_{j-1}\mon w_j} \,\a_j\r)  
\end{eqnarray*}
in $\Fs(x\sox\ul w\sox y)$ with $\eta_{j-1,j} := \eta_{w_0}\circ\cdots \circ\eta_{w_{j-2}}\circ \eta_{w_{j+1}}\circ \cdots \circ\eta_{w_m}$. In the case $w_0 = w_1 = \cdots = w_{m+1}=x$ and $\Fs$ is injective, this reduces in $\Fs(x\sox x)$ to
\begin{equation}\label{log a1am 2}
\log_{x\mon x} (\a_{m+1}\a_m\cdots\a_1)  \=  \sum_{j=1}^{m+1} \log_{x\mon x} \,\a_j.
\end{equation}
\end{prop}

Proof:\;  For 1. one has 
\begin{eqnarray*}
\log_{x\mon x \mon x} (x \stackrel{p}{\to} x\stackrel{p}{\to} x)  & = &\eta_\lx \log_{x\mon x} (x \stackrel{p}{\to} x) + \eta_\rx \log_{x\mon x} (x\stackrel{p}{\to} x)  \\[1mm] 
& \stackrel{p\circ p \,=\,p}{=} & \eta_\lx \log_{x\mon x} (x \stackrel{p\circ p}{\to} x) + \eta_\rx\log_{x\mon x} (x\stackrel{p\circ p}{\to} x)  \\[1mm] 
& \stackrel{\eqref{log compatible 2}, \eqref{log compatible 3}}{\=} & \log_{x\mon x \mon x} (x \stackrel{p}{\to} x\stackrel{p}{\to} x) + \log_{x\mon x \mon x} (x \stackrel{p}{\to} x\stackrel{p}{\to} x).
\end{eqnarray*}
Hence 
$ 0 \=  \log_{x\mon x \mon x} (x \stackrel{p}{\to} x\stackrel{p}{\to} x) \stackrel{\eqref{log compatible 2}}{\=} \eta_\lx (\log_{x\mon x} \,p\circ p)  = \eta_\lx (\log_{x\mon x} \, p).$
The other statements follow similarly.
\begin{flushright}
$\Box$
\end{flushright}

{\bf Comments:} 
If the pretracial monoidal product representation $\Fs: \Cb^* \to  \Ring_{\mbox{{\tiny {\bf Add}}}}$ is endowed with a trace $\t$ then the $\t$-character of the log-functor defines a {\it log-determinant functor representation} of $\Cb$, mapping each $ w\in\ob(\Cb)$ to  
$\emor_\Ab(1)$ and $\a\in\mor_{\ul z}(x,y)$ to the character  
\begin{equation*}\label{log det tau}
\wt\t(\log\, \a)  := \wt\t_{x\mon\ul z\mon y}(\log_{x\mon\ul z\mon y}\, \a) \in \emor_\Ab(1),
\end{equation*}
of $\log_{x\mon\ul z\mon y}\,\a \in \Fs(x\mon\ul z\mon y)/[\Fs(x\mon\ul z\mon y), \Fs(x\mon\ul z\mon y)]$. 
The character of $\a\in\mor_z(x,y)\in\Nn_p\Cb$ is invariantly defined: in $\mor_\Ab(1,1)$ 
\begin{equation}\label{trace log unique}
\wt\t_{x\mon z \mon y}\l(\log_{x\mon z\mon y} \,\a\r)  = \wt\t_{x\mon y } \l(\log_{x\mon y} \,\a\r).
\end{equation}
Likewise, for $\d\in\mor(x,y)$
$\wt\t_{x\mon z \mon y}\l(\eta_z (\log_{x\mon y} \,\d)\r)  = \wt\t_{x\mon y } \l(\log_{x\mon y} \,\d\r)$, and more generally  with  $\ul z= (z_1, \ldots, z_r)$, $x= x_1\sox \cdots \sox x_n$ one has
\begin{equation}\label{trace log unique 2}
\wt\t_{x_{\ul z}}\l(\log_{x_{\ul z}} \,\nu\r)  =   \wt\t_x\l(\log_x \, \nu\r).
\end{equation}
Indeed, for   $w \in \ob(\Cb)$ one has $\log_{x_{\ul w}} (\nu)  -   \eta_{\ul w}(\log_{x_1\mon \cdots \mon x_n} \,\nu) \ \in \ [\Fs(x_{\ul w}), \Fs(x_{\ul w})]$ by \eqref{log collapsing two simplices} whilst   $[\Fs(w), \Fs(w)] \<\Ker (\t_w)$. Hence   \eqref{two traces 2} yields the conclusion.

Here, \eqref{trace log unique} is shorthand for
$\wt\t_{x\mon z \mon y}\l(\log_{x\mon z\mon y}  (x \stackrel{\a}{\to} z \stackrel{\b}{\to} y)\r)  = \wt\t_{x\mon y }\l(\log_{x\mon y} (x \stackrel{\b\circ\a}{\to} y)\r)$,  \ or $\wt\t_{x\mon z \mon y}(\log_{x\mon z\mon y}  \,\b\a)  = \wt\t_{x\mon y }(\log_{x\mon y} d_1(\b\a ))$. By the above, \emph{the} logarithmic character 
 \eqref{log det tau}, 
of a morphism $\a \in \mor_\Cb(x,y)$ is independent of where it is computed.

For $\a\in \mor(x,z)$ and $\b\in \mor(z,y)$ one has
$
\wt\t\l(\log\,\b\a\r) = \wt\t\l(\log \,\a\r) + \wt\t\l(\log \,\b\r)$ in $\mor_\Ab(1,1)$.

The space $\Log(\Cb, \Fs)$  of logarithms on $\Cb$ with respect to a fixed monoidal product representation $\Fs$ is an abelian group
$\log_1, \log_2 \in  \Log(\Cb, \Fs)   \To  \log_1 + \log_2  \in  \Log(\Cb, \Fs)$
with respect to the additive structure of the category $\Ab$, 
as is the space $\Log^\chi(\Cb)$   of logarithmic characters $\t(\log \a)$ independently of a particular $\Fs$

If $\Cb$ is an additive category then $\t\circ \log$ is a log-representation from the maximal sub groupoid of $\Cb$, whose morphisms are the isomorphisms of $\Cb$, to the isomorphism torsion group $K_1^{{\rm {\tiny iso}}}(\Cb)$ of \cite{Ra}.   

By statement 5 (and 6) of \propref{log properties} it is enough to require log-additivity on 1-simplicies to infer it on $p$-simplices in $\Nn\Cb$. On the other hand, as far as computing log-determinant characters is concerned,  log-additivity \eqref{log equality 2} can be formulated more generally as the existence of  $\ul w_0, \ul w_1, \ul w_2 \in \ob(\Cb)$ such that $ \wt\eta_{\ul w_0}(\log_{x\mon z} \,\a),  \wt\eta_{\ul w_1}(\log_{z\mon y}\, \b), \wt\eta_{\ul w_2}(\log_{x\mon y} (x \stackrel{\b\circ\a}{\to} y))$ are all in the same $\Fs(v)$ with, in $ \Fs(v)/[ \Fs(v),  \Fs(v)],$   
\begin{equation}\label{log equality p} 
\wt\eta_{\ul w_1}(\log_{x\mon y} (x \stackrel{\b\circ\a}{\to} y) )  = \wt\eta_{\ul w_2}(\log_{x\mon z} \,\a) + \wt\eta _{\ul w_0}(\log_{z\mon y}\, \b).
\end{equation} 

Despite \lemref{1-simplices rule ok}, it can be natural to define simplicial logarithms directly on strata $\mor_{\ul z}(x,y)$ in $p$-simplices with $p>1$. In particular,  this allows a log-functor to be extended to  $\d\in \mor_\Cb(1,1)=\emor_\Cb(1)$ factorisable as  $\d = \b\a$ for  $\a\in \mor_\Cb(1,z)$ and $\b\in \mor_\Cb(z,1)$ with $z\neq 1\in\ob(\Cb)$ (this is always the case on $\cob_n$). Choosing such a factorisation, define
\begin{equation}\label{logs on end(1)}
\log_z \,\d := \log_z (1\stackrel{\a}{\to} z \stackrel{\b}{\to} 1) \ \in \Fs(z)/[ \Fs(z),  \Fs(z)].
\end{equation}
Here, we use $\log_z := \log_{1\mon z\mon 1}$ and $\Fs(1\sox z\sox 1) = \Fs(z)$, as $\Fs$ is exact and $\log$ is strict, which depends on $\d$ and $z$ but   is independent of the particular choice of $\a, \b$. In the presence of a trace one then further has
\begin{equation}\label{log1}
\log_1:\emor_\Cb(1) \to  \emor_\Ab(1),  \hf  \log_1 \,\d := \wt\t(\log_z (1\stackrel{\a}{\to} z \stackrel{\b}{\to} 1)),
\end{equation}
independently of the particular choice of $\a, \b$ and of $z$ and by log-additivity
\begin{equation}\label{log2}
  \log_1 \,\d := \wt\t(\log_z \,\a) + \wt\t(\log_z \,\b)
\end{equation}
as a particular case of the additivity of log-characters.

\section{Log-structures on bordism categories}

There are a number of bordism categories with natural logarithmic functors.  Bordism classes  will be denoted  $\ol W \in\mor_{{\bf Bord_n}}(M_0,M_1)$, while $W=(W,\kappa_{\mbox{{\tiny $\pd W$}}})\in \ol W $ will indicate a smooth representative of the class. Thus, $W$ is an oriented smooth compact manifold of dimension $n+1$ whose boundary $\pd W\in\ob(\cob_n)$ is endowed with an orientation preserving diffeomorphism 
$\kappa_{\mbox{{\tiny $\pd W$}}}: \pd W \to M_0^- \sqcup M_1$, the superscript indicating the reverse orientation on $M_0$. $\ol W= \ol{(W,\kappa_{\mbox{{\tiny $\pd W$}}})} $ denotes the equivalence class relative to oriented diffeomorphism. 

Let
$\Fs: \cob^*_n\to \Ring_{\mbox{{\tiny {\bf Add}}}}$
be an unoriented pretracial monoidal product representation.  {\em Unoriented} is the assumption that
$
\Fs(M^{(-)})  = \Fs(M),
$
where $M^{(-)}$ denotes $M$ with one or more of its connected components with orientation reverse.
A log-TQFT on $\cob_n$ relative to $\Fs$ is a log-additive presimplicial map 
\begin{equation*}
\log: \Nn\cob_n\to \Fs(\cob_n^*)/[\Fs(\cob_n^*), \Fs(\cob_n^*)],  
\end{equation*}
defining  for each $p$-simplex 
$M_0\stackrel{\mbox{{\tiny $\ol W_0$}}}{\to} M_1 \stackrel{\mbox{{\tiny $\ol W_1$}}}{\to} M_2 \to \cdots \to M_{p-1} \stackrel{\mbox{{\tiny $\ol W_{p-1}$}}}{\to} M_p \in \Nn_p\cob_n$ of bordisms between compact boundaryless manifolds $M_j$, a logarithm 
\begin{equation}\label{log on bordisms}
\log_M(M_0\stackrel{\mbox{{\tiny $\ol W_0$}}}{\to}  M_1 \stackrel{\mbox{{\tiny $\ol W_1$}}}{\to}  M_2 \to \cdots \to M_{p-1} \stackrel{\mbox{{\tiny $\ol W_{p-1}$}}}{\to}  M_p)  \ \in \ \Fs_{\mbox{{\tiny $\Pi$}}}( M) := \Fs( M)/[\Fs( M), \Fs( M)],
\end{equation}
where $ M =  M_0 \sqcup M_1 \sqcup\cdots\sqcup M_p$, with  
\begin{equation}\label{logtqft simplicial}
\log_{\mbox{{\tiny $M_0\sqcup M_1 \sqcup M_2$}}}(M_0\stackrel{\mbox{{\tiny $\ol W_0$}}}{\to} M_1\stackrel{\mbox{{\tiny $\ol W_1$}}}{\to} M_2 )
 = \wt\eta_{\mbox{{\tiny $M_1$}}}
\log_{\mbox{{\tiny $M_0\sqcup M_2$}}}(M_0\stackrel{\mbox{{\tiny $\ol W_0 \cup \ol W_1$}}}{\too} M_2 ),
\end{equation}
where  $\ol W_0 \cup \ol W_1$ is the composed bordism joined along $M_1$, and, on 1-simplices,
\begin{equation}\label{logtqft add}
\wt\eta_{\mbox{{\tiny $M_1$}}}\log_{\mbox{{\tiny $M_0\sqcup M_2$}}}(\ol W_0 \cup \ol W_1) = \wt\eta_{\mbox{{\tiny $M_2$}}}\log_{\mbox{{\tiny $M_0\sqcup M_1$}}}(\ol W_0) +  \wt\eta_{\mbox{{\tiny $M_0$}}}\log_{\mbox{{\tiny $M_1\sqcup M_2$}}}(\ol W_1)
\end{equation}
in $\Fs_{\mbox{{\tiny $\Pi$}}}(M_0\sqcup M_1 \sqcup M_2).$
Though $\Fs$ is unoriented, the logarithms $\log_{\mbox{{\tiny $M$}}}(\ol W) $ {\em will} in  general depend on the orientation of the bordisms $\ol W$.
 The $M_j$ need not be connected. On the other hand, writing $M_j = N_0 \sqcup \cdots \sqcup N_k$   is reflected functoriality  in a canonical isomorphism 
$\Fs(M_j) \cong  \Fs(N_0\sqcup \cdots \sqcup N_k).$
A permutation of the ordering $N_{\s(0)} \sqcup \cdots \sqcup N_{\s(k)}$ yields (in accordance with \eqref{F symmetric 1}) a compatible canonical isomorphism
$
\mu_\s : \Fs(N_0 \sqcup \cdots \sqcup N_k)  \stackrel{\cong}{\to} \Fs(N_{\s(0)} \sqcup \cdots \sqcup N_{\s(k)}).
$ In \eqref{log on bordisms} there is no ambiguity because $M$ is defined to be the given disjoint union in the order specified by the $p$-simplex. 

The p-simplices of $\Nn\cob_n$ may be viewed as  bordisms which retain  data of how they were formed by gluing other bordisms. Boundaryless bordisms $\ol W\in\mor_{{\bf Bord_n}}(\mt, \mt)$ need separate consideration: we are instructed by \eqref{logs on end(1)} to  view $\ol W$ as a composed bordism $\mt \stackrel{\mbox{{\tiny $\ol W_0$}}}{\too} M \stackrel{\mbox{{\tiny $\ol W_1$}}}{\too} \mt $ relative to codimension 1 embedded submanifold $M \hookrightarrow W$ and set
$$
\log_{\mbox{{\tiny $M$}}} \, \ol W :=\log_{\mbox{{\tiny $M$}}}(\mt \stackrel{\mbox{{\tiny $\ol W_0$}}}{\too} M \stackrel{\mbox{{\tiny $\ol W_1$}}}{\too} \mt )  \hf\in\ \Fs(M)/[\Fs(M), \Fs(M)].
$$
Log-additivity then gives 
$
\log_{\mbox{{\tiny $M$}}} \,\ol W = \log_{\mbox{{\tiny $M$}}}(\mt\stackrel{\mbox{{\tiny $\ol W_0$}}}{\to} M ) + \log_{\mbox{{\tiny $M$}}}(M\stackrel{\mbox{{\tiny $\ol W_1$}}}{\to} \mt ) \in \Fs(M)/[\Fs(M), \Fs(M)],
$
and if  tracial with character
$
\t(\log \ol W) = \t_{\mbox{{\tiny $M$}}}( \log \ol W_0) + \t_{\mbox{{\tiny $M$}}}( \log \ol W_1)  \in \ \emor_\Ab(1)
$
depending only on $\ol W$, not on its factorisation as $\ol W_0 \cup_{\mbox{{\tiny $M$}}} \ol W_1$.

\begin{lem}\label{log cylinder}  
Let $C_M\in\mor_{{\bf Bord_n}}(M,M)$ be the bordism class of $[0,1]\x M$. Then  
\begin{equation*}
\wt\eta_{\mbox{{\tiny $M$}}}\log_{\mbox{{\tiny $M\sqcup M$}}} (C_M) = 0,
\end{equation*}
in $\Fs_{\mbox{{\tiny $\Pi$}}}(M\sqcup M\sqcup M)$ and $\log_{\mbox{{\tiny $M\sqcup M$}}} (C_M) = 0 \in \Fs_{\mbox{{\tiny $\Pi$}}}(M\sqcup M)$ if $\Fs$ is injective. For $\ol W\in \mor(M_0, M_1)$ 
\begin{equation}\label{logtqft}
\wt \eta_{\mbox{{\tiny $\sqcup N$}}}\, \log_{\mbox{{\tiny $M_0 \sqcup M_1$}}}(M\stackrel{\mbox{{\tiny $\ol W$}}}{\to}N) =   \log_{\mbox{{\tiny $M_0 \sqcup M_1 \sqcup N$}}}(M\stackrel{\mbox{{\tiny $\ol W$}}}{\to}N\stackrel{\mbox{{\tiny $C_{\mbox{{\tiny $N$}}}$}}}{\to}N)
\end{equation}
in $\Fs_{\mbox{{\tiny $\Pi$}}}(M_0 \sqcup M_1\sqcup N)$.
\end{lem}
Proof:\;   Restatements of \propref{log properties} (1) and (2) to $\cob_n$. $\hfill \Box\\[1mm]$

A log-TQFT yields a TQFT,  in the following sense: 
\begin{lem}\label{log tqft => tqft}  
A log-TQFT, defined by $\log: \Nn\cob_n\to \Ring_{\mbox{{\tiny {\bf Add}}}}$ relative to a tracial $\Fs: \cob_n^* \to (\Fs(\cob_n^*), \t)$ defines a monoid $(\mor_\Ab(1,1), +)$-valued symmetric monoidal functor $Z_{{\rm log},\t,\e}: \cob_n\to \mor_\Ab(1,1)$  by setting $Z_{{\rm log},\t,\e}(M)=\mor_\Ab(1,1)$ and  $Z_{{\rm log},\t,\e}(\ol W) = \t(\log \ol W).$
\end{lem}

Conversely, log-TQFTs may arise from TQFTs, but we know of this  in essentially  trivial cases only. Non-trivial log-TQFTS are not hard to find, however. 

Let $\cob^*_n$  be the subcategory of  $\cob_n$ whose morphisms are the coherence and permutation bordisms. Define a monoidal product  representation 
$
\Fs_{\mbox{{\tiny $-\oo$}}}: \cob_n^*\to \alg_{\texttt{F}}$ by setting
$\Fs_{\mbox{{\tiny $-\oo$}}}( M) := \Psi^{-\oo}(  M) :=  \Psi^{-\oo}(  M, \wedge T^* M)$ 
to be the algebra of smoothing operators  on  the de Rham complex $\O(M)$ with the coherence bordisms of the monoidal product $\sqcup$ mapped to the identity operator. An element $T\in \Fs_{\mbox{{\tiny $-\oo$}}}( M)$ is  specified by a Schwartz kernel
\begin{equation}\label{smoothing kij}
 k_{\mbox{{\tiny $M$}}} \in C^\oo(M \x M, ((\wedge T^*M)^* \ox |\Lambda|^{\mbox{{\tiny $\frac{1}{2}$}}}_{\mbox{{\tiny $M$}}}) \boxtimes (\wedge T^*M\ox|\Lambda|^{\mbox{{\tiny $\frac{1}{2}$}}}_{\mbox{{\tiny $M$}}}))
\end{equation}
taking values in form valued half-densities
 
If $M$ is disconnected and is written as a disjoint union $M =M_1 \sqcup \cdots \sqcup M_m$  of $M_j\in \ob(\cob_n)$, then 
$\O(M) = \O(M_1)  \oplus \cdots \oplus  \O(M_m)$ with respect to which $T\in \Fs_{\mbox{{\tiny $-\oo$}}}( M)$ is an $n\x n$ block matrix $(T_{i,j})$ of smoothing operators $T_{i,j}\in \Psi^{-\oo}(M_j , M_i)$  specified  by Schwartz kernels 
\begin{equation}\label{smoothing kij}
 k_{i,j} \in C^\oo(M_i \x M_j, ((\wedge T^*M_i)^* \ox |\Lambda|^{\mbox{{\tiny $\frac{1}{2}$}}}_{\mbox{{\tiny $M_i$}}}) \boxtimes (\wedge T^*M_j \ox|\Lambda|^{\mbox{{\tiny $\frac{1}{2}$}}}_{\mbox{{\tiny $M_j$}}}))
\end{equation}
whose rows and columns are permuted by $\mu_\s( M)$ relative to a reordering $\s$ of the $M_j$. 

With $i: M :=M_1 \sqcup \cdots \sqcup M_m \hookrightarrow M_{\mbox{{\tiny $N$}}} := M_1 \sqcup \cdots \sqcup N \sqcup \cdots\sqcup M_m$, the insertion maps are the canonical inclusions
\begin{equation}\label{etaN}
\eta_{\mbox{{\tiny $N$}}} : \Fs_{\mbox{{\tiny $-\oo$}}}(  M)   \hookrightarrow \Fs_{\mbox{{\tiny $-\oo$}}}( M_{\mbox{{\tiny $N$}}}), \hf \eta_{\mbox{{\tiny $N$}}}(T) = i_{\mbox{{\tiny $N$}}} \circ T\circ i_{\mbox{{\tiny $N$}}}^*.
\end{equation}
$\Fs_{\mbox{{\tiny $-\oo$}}}$ is pretracial, though not injective, and we may form the pushed-down insertion maps 
 \begin{equation}\label{pushed-down insertion maps}
\wt\eta_{\mbox{{\tiny $N$}}} = \wt\eta_{\mbox{{\tiny $N$}}}(M) : \frac{\Fs_{\mbox{{\tiny $-\oo$}}}( M)}{[\Fs_{\mbox{{\tiny $-\oo$}}}( M), \Fs_{\mbox{{\tiny $-\oo$}}}( M)]}   \to \frac{\Fs_{\mbox{{\tiny $-\oo$}}}( M_{\mbox{{\tiny $N$}}})}{[\Fs_{\mbox{{\tiny $-\oo$}}}( M_{\mbox{{\tiny $N$}}}), \Fs_{\mbox{{\tiny $-\oo$}}}( M_{\mbox{{\tiny $N$}}})]}.
\end{equation}

\begin{lem}\label{Tr unique} 
The linear map
 \begin{equation}\label{cobordism mpr trace}
\Tr_{\mbox{{\tiny $ M$}}} : \Fs_{\mbox{{\tiny $-\oo$}}}( M) \to \C, \hf \Tr_{\mbox{{\tiny $ M$}}} (T) :=  \sum_{j=1}^m \Tr_{\mbox{{\tiny $M_j$}}} (T_{j,j}) :=  \sum_{j=1}^m\int_{M_j} \tr(k_{j,j}(m,m)),
\end{equation}
is a trace and, up to a multiplication by a constant, is the unique trace on $\Fs_{\mbox{{\tiny $-\oo$}}}( \cob_n^*)$. 
The quotients $\frac{\Fs_{\mbox{{\tiny $-\oo$}}}( M)}{[\Fs_{\mbox{{\tiny $-\oo$}}}( M), \Fs_{\mbox{{\tiny $-\oo$}}}( M)]}$ are complex lines and the trace  defines and is defined by a linear isomorphism 
\begin{equation}\label{trace reduced} 
\widetilde \Tr_{\mbox{{\tiny $M$}}} : \frac{\Fs_{\mbox{{\tiny $-\oo$}}}( M)}{[\Fs_{\mbox{{\tiny $-\oo$}}}( M), \Fs_{\mbox{{\tiny $-\oo$}}}( M)]} \ \stackrel{\cong}{\to} \ \C
\end{equation}
with 
\begin{equation}\label{commuting trace reduced} 
\Tr_{\mbox{{\tiny $ M$}}} =  \wt\Tr_{\mbox{{\tiny $ M$}}}\circ \pi_{\mbox{{\tiny $ M$}}}.
\end{equation}
One has 
\begin{equation}\label{tracial Finfty} 
\Tr_{\mbox{{\tiny $ M$}}} =  \Tr_{\mbox{{\tiny $ M_{\mbox{{\tiny $N$}}}$}}}\circ \eta_{\mbox{{\tiny $N$}}} \hf \oon \ \Fs_{\mbox{{\tiny $-\oo$}}}( M),\hf\hf\hf\hf
\end{equation}
\begin{equation}\label{tilde tracial Finfty} 
\ \hf\hf\hf\wt\Tr_{\mbox{{\tiny $ M$}}} =  \wt\Tr_{\mbox{{\tiny $ M_{\mbox{{\tiny $N$}}}$}}}\circ \wt\eta_{\mbox{{\tiny $N$}}}\hf \ \oon \ \Fs_{\mbox{{\tiny $-\oo$}}}( M)/[\Fs_{\mbox{{\tiny $-\oo$}}}( M), \Fs_{\mbox{{\tiny $-\oo$}}}( M)].
\end{equation}
\end{lem} 

We omit the straightforward proof. 
 
The pushed-down insertion map $\wt\eta_{\mbox{{\tiny $N$}}} (M)$ in \eqref{pushed-down insertion maps}  is hence a linear isomorphism of complex lines, and fits into the commutative diagram \eqref{trace diagram} which, here, is 

\begin{equation}\label{trace diagram Tr infty} 
\begin{array}{rcccl}
  \Fs_{\mbox{{\tiny $-\oo$}}}(M) &  & \stackrel{ \eta_{\mbox{{\tiny $N$}}}(M)}{\longrightarrow} &  & \  \Fs_{\mbox{{\tiny $-\oo$}}}(M_{\mbox{{\tiny $N$}}})   \\
   & \stackrel{\ \ \ \ \Tr_{\mbox{{\tiny $ M$}}}}{\searrow} & \hf \hf\hf\hf &  \stackrel{\Tr_{\mbox{{\tiny $M_{\mbox{{\tiny $N$}}}$}}}}{\hf\ \ \swarrow}  \hskip 5mm &    \\
    &  & &  &   \\
  \downarrow \pi_{\mbox{{\tiny $ M$}}}   &  & \C &  & \hf\downarrow \pi_{\mbox{{\tiny $ M_{\mbox{{\tiny $N$}}}$}}}  \\
      &  & &  &  \\
  & \stackrel{\wt\Tr_{\mbox{{\tiny $ M$}}}}{\ \ \ \hf\nearrow}   \hf\hf & & \stackrel{\ \ \ \wt\Tr_{\mbox{{\tiny $ M_{\mbox{{\tiny $N$}}}$}}}}{\nwarrow}  &   \\   \frac{ \Fs_{\mbox{{\tiny $-\oo$}}}(M)}{[ \Fs_{\mbox{{\tiny $-\oo$}}}(M), \,  \Fs_{\mbox{{\tiny $-\oo$}}}(M)]}  &  & \stackrel{\wt\eta_{\mbox{{\tiny $N$}}} (M)\ \cong}{\longrightarrow} &   & \frac{ \Fs_{\mbox{{\tiny $-\oo$}}}(M_{\mbox{{\tiny $N$}}})}{[ \Fs_{\mbox{{\tiny $-\oo$}}}(M_{\mbox{{\tiny $N$}}}),  \, \Fs_{\mbox{{\tiny $-\oo$}}}(M_{\mbox{{\tiny $N$}}})]}  
\end{array},
\end{equation}
and one has
$
\wt\eta_{\mbox{{\tiny $N$}}} (M)= \wt\Tr_{\mbox{{\tiny $ M_{\mbox{{\tiny $N$}}}$}}}\ii\circ \wt\Tr_{\mbox{{\tiny $ M$}}}.$ 
Likewise, in view of the isomorphism \eqref{trace reduced}, $\pi_{\mbox{{\tiny $ M$}}}(A)$ may be characterised as {\it the abstract trace} of $A\in \Fs_{\mbox{{\tiny $-\oo$}}}(M)$, one has
$
\pi_{\mbox{{\tiny $ M$}}} = \wt\Tr_{\mbox{{\tiny $ M$}}}\ii\circ\Tr_{\mbox{{\tiny $ M$}}}.$

The classical trace hence refines $\Fs_{\mbox{{\tiny $-\oo$}}}$ to a tracial monoidal product representation $(\Fs_{\mbox{{\tiny $-\oo$}}}, \Tr)$.   There is, on the other hand,  the `larger' monidal product representation 
\begin{equation}\label{Z tmpr}
\Fs_{\mbox{{\tiny $\Z, -\oo$}}}: \cob_n^*\to \alg_{\texttt{F}}, \hf  M\mto \Fs_{\mbox{{\tiny $\Z, -\oo$}}}( M)
\end{equation}
with    $\Fs_{\mbox{{\tiny $\Z, -\oo$}}}( M) $ the algebra of continuous operators on  $\O( M)$ defined by Schwartz kernels which are smoothing off the `matrix diagonal' and pseudodifferential along it, in the following sense. Let $M_1, \ldots, M_m$ be the connected components of $M$ and let $k_{i,j}$ be the restriction to $M_i\x M_j$ of the distributional kernel of $T\in \Fs_{\mbox{{\tiny $\Z, -\oo$}}}( M)$.  Then $k_{i,j}$ is required to be a smoothing kernel \eqref{smoothing kij} if $i\neq j$, while if $i=j$ it may, more generally, be an integer order  pseudodifferential operator  ($\pdo$) kernel 
$k_{j,j} \in \Dd\pr(M_j \x M_j, ((\wedge T^*M_j)^* \ox |\Lambda|^{\mbox{{\tiny $\frac{1}{2}$}}}_{\mbox{{\tiny $M_j$}}}) \boxtimes (\wedge T^*M_j \ox|\Lambda|^{\mbox{{\tiny $\frac{1}{2}$}}}_{\mbox{{\tiny $M_j$}}}))$
 in the space of conormal distributions on form-valued half-densities. Thus, there is an atlas of $M_j \x M_j$ in which $k_{j,j}$ can be written in each localisation as an oscillatory integral  
\begin{equation}\label{osc int}
k_{j,j}(x,y) =  \int_{\R^n} e^{i\xi.(x-y)}\, {\bf b} ^{[j]}(x, y, \xi) \ \db \xi \, |dx|^{\mbox{{\tiny $\frac{1}{2}$}}} |dy|^{\mbox{{\tiny $\frac{1}{2}$}}}
\end{equation}
of a symbol (amplitude) ${\bf b}^{[j]}(x,y,\xi)$ of order $p_j\in\Z\cup\{-\oo\}$ (depending on the trivialisation).   $\Fs_{\mbox{{\tiny $\Z, -\oo$}}}( M)$ is filtered by the  subspaces 
$ 
\Fs_{\mbox{{\tiny $p, -\oo$}}}( M)= \Psi^{\mbox{{\tiny $p, -\oo$}}}( M)
$
of operators with classical $\pdos$ on the diagonal up to order $p\in \Z$. If $M =M_1 \sqcup \cdots \sqcup M_m$ then $\Fs_{\mbox{{\tiny $\Z, -\oo$}}}( M)$ is identified with the matrix algebra  $(T_{i,j})$ of  operators $T_{i,j}$ with smoothing kernels off the matrix diagonal and with integer order $\pdo$ oscillatory kernel \eqref{osc int} if $i=j$.   

$\Fs_{\mbox{{\tiny $\Z, -\oo$}}}$ is pretracial with quotient functor $\rho_{\mbox{{\tiny $ M$}}} : \Fs_{\mbox{{\tiny $\Z, -\oo$}}}( M) \to  \Fs_{\mbox{{\tiny $\Z, -\oo$}}}( M)/[\Fs_{\mbox{{\tiny $\Z, -\oo$}}}( M), \Fs_{\mbox{{\tiny $-\oo$}}}( M)]$. It has a trace structure complementary to the classical trace and not quite unique:
\begin{lem}\label{res m-unique} Let $M_j$  be the  connected components of $M$. Then the linear space of traces  on $\Fs_{\mbox{{\tiny $\Z,  -\oo$}}}( M)$ has (complex)  dimension $m$: on $\Fs_{\mbox{{\tiny $\Z, -\oo$}}}( M)$ each ${\bf c} = (c_1, \ldots , c_m)\in \C^m$ parametrises the linear sum of residue traces 
\begin{equation}\label{res m-trace}
 \res^{{\bf c}}_{\mbox{{\tiny $ M$}}} (B) = \sum_{j=1}^m c_j \res_{\mbox{{\tiny $M_j$}}}(B_{jj}) :=  \sum_{j=1}^m c_j \int_{S^*M_j} b^{[j]}_{-n}(x,\eta) \, \db_S\eta\, |dx|.
\end{equation}
Each such trace  defines and is defined by a linear homomorphism 
\begin{equation}\label{res trace reduced} 
\wt\res^{{\bf c}}_{\mbox{{\tiny $ M$}}} : \frac{\Fs_{\mbox{{\tiny $\Z, -\oo$}}}( M)}{[\Fs_{\mbox{{\tiny $\Z, -\oo$}}}( M), \Fs_{\mbox{{\tiny $\Z, -\oo$}}}( M)]} \ \stackrel{\cong}{\to} \ \C \ \ 
\wwith \ \ 
\res^{{\bf c}}_{\mbox{{\tiny $ M$}}}=  \wt\res^{{\bf c}}_{\mbox{{\tiny $ M$}}}\circ \rho_{\mbox{{\tiny $ M$}}}.
\end{equation}
\end{lem}

These structures behave well with respect to diffeomorphisms:

\begin{lem}\label{phi sharp}
Let $\Fs: \cob_n^*\to \alg_{\texttt{F}}, \ M\mto (\Fs(M), \t_{\mbox{{\tiny $M$}}}),$ be either one of the tracial monoidal product representations $ (\Fs_{\mbox{{\tiny $-\oo$}}}, \Tr)$ or $ (\Fs_{\mbox{{\tiny $\Z, -\oo$}}}, \res)$. Let $M^{(-)}$ be $M$ with one or more of its connected components with orientation reversed.  Then
$
\Fs\l(M^{(-)}\r) \  =  \ \Fs(M). $
A diffeomorphism $\phi: M\to N$ between $M, N \in \ob(\cob_n)$ induces a canonical continuous isomorphism of algebras  
$ 
\phi_\sharp : \Fs( M) \  \to  \ \Fs( N), 
$
preserving the filtration by $\psdo$ order, and pushes-down to a continuous linear map $\wt\phi_{\mbox{{\tiny $M, N$}}} : \Fs(M) / [\Fs(M), \Fs(M)] \to \Fs(N) / [\Fs(N), \Fs(N)].$

Trace invariance: there is a commutative diagram
\begin{equation}\label{trace diagram theta mn} 
\begin{array}{rcccl}
  \Fs(M) \ &  & \stackrel{\phi_\sharp}{\longrightarrow} &  & \  \Fs(N)   \\
   & \stackrel{\t_{\mbox{{\tiny $M$}}}}{\searrow} & \hf \hf\hf\hf &  \stackrel{\t_{\mbox{{\tiny $N$}}}}{\hf\swarrow}  \hskip 5mm &    \\
    &  & &  &   \\
  \downarrow\pi_{\mbox{{\tiny $M$}}}  &  & \C &  & \hf\downarrow\pi_{\mbox{{\tiny $N$}}} \\
    &  & &  &  \\
  & \stackrel{\tilde\t_{\mbox{{\tiny $M$}}}}{\hf\nearrow} \hskip 5mm&  & \stackrel{\tilde\t_{\mbox{{\tiny $N$}}}}{\nwarrow} &    \\
 \frac{\Fs(M)}{[\Fs(M), \Fs(M)]}  &  & \stackrel{\wt\phi_\sharp }{\longrightarrow} &  & \frac{\Fs(N)}{[\Fs(N), \Fs(N)]}  
\end{array}.
\end{equation}

For $ (\Fs_{\mbox{{\tiny $-\oo$}}}, \Tr)$ the map $\wt\phi_\sharp$ is independent of the choice of $\phi$: if $M$ and $N$ are diffeomorphic there is a canonical linear isomorphism of complex lines:
\begin{equation}\label{phi sharp 2} 
\vartheta_{\mbox{{\tiny $M, N$}}} : \frac{\Fs_{\mbox{{\tiny $-\oo$}}}(M)}{[\Fs_{\mbox{{\tiny $-\oo$}}}(M), \Fs_{\mbox{{\tiny $-\oo$}}}(M)]}   \ \to \ \frac{\Fs_{\mbox{{\tiny $-\oo$}}}(N)}{[\Fs_{\mbox{{\tiny $-\oo$}}}(N), \Fs_{\mbox{{\tiny $-\oo$}}}(N)]}.
\end{equation}
\end{lem}

This is readily checked; thus, the diffeomorphism  $\phi$ induces a bundle isomorphism $\wedge TN^*\to \wedge TM^*$  and hence a continuous linear pull-back isomorphism $\phi_* : \O(N) \stackrel{\cong}{\to} \O(M)$, with respect to which $\phi_\sharp (T) :=  \phi_*\ii\circ T \circ \phi_*$ is an algebra isomorphism  defining  an abelian groupisomorphism 
$
 [\Fs(M), \Fs(M)]   \ \stackrel{\cong}{\to} \ [\Fs(N), \Fs(N)]$.
which with \eqref{phi sharp} gives \eqref{phi sharp 2}. For the diagram,one uses the universality property of traces and  Lidskii's theorem.

\subsubsection{The topological signature}

 For a compact oriented manifold  $W$  of dimension $4k$ with boundary $\pd W$,  the topological signature $\sgn(W)$ of $W$, defined to be the signature of the quadratic form 
\begin{equation}\label{signature}
\wh H^{2k}(W) \x \wh H^{2k}(W) \to \R, \hf (\xi, \xi\pr) \mto <\xi\cup\xi\pr, [W]>,
\end{equation}
with $\wh H^{2k}(W)$ the image of the inclusion  $ H^{2k}(W,\pd W)  \to H^{2k}(W)$
This  arises as a character of a logarithmic representation on bordisms as follows. 

On a smooth representative $W\in \ol W$ of a bordism class $\ol W \in \mor_{\mbox{{\tiny $ {\rm Bord}_{4k}$}}}(M_0, M_1)$, a choice of Riemannian metric $g_W$ is made which in a collar neighbourhood $U_j$ of each boundary component $\pd W_j$ is a product metric $g_{U_j}  = du_j^2 + g_{\partial W_j}$ with $u_j$ a choice of normal coordinate in $(-1,0]$   if $\pd W_j$ is a component of $M_0^-$ and in  $[0,1)$   if $\pd W_j$ is a component of $M_1$; all logarithms will be independent of the choice of $g_W$ and the choice of representative $W$. Associated to  $g_W$ is a Hodge star isomorphism 
$*: \O^p(W) \to \O^{4k-p}(W)$ and a signature operator $$\eth^{\mbox{{\tiny$W$}}} = d + d^*: \O^+(W) \to \O^-(W)$$ between
the eigenspaces $\O^\pm(W)$ of the involution $i^{p(p-1)} *$ on the de Rham complex. 

Recall from \cite{APS1},  since $W$ is isometric to a product near each boundary component $\pd W_j$ the operator $\eth^{\mbox{{\tiny$W$}}}$ acts along tangential boundary directions by a self-adjoint  signature operator $B_j$ on the de Rham algebra $\Omega(\pd W_j)$, equal  to 
$B^{2p}_j := (-1)^{k+p+1}(*d_j - d_j*)$  on $\Omega^{2p}(\pd W_j)$ and to $B^{2p-1}_j := (-1)^{k+p} (*d_j + d_j*)$ on $\Omega^{2p-1}(\pd W_j)$.   Let $B_j^{ev} = \bigoplus_p B^{2p}_j$, $B_j^{odd} = \bigoplus_p B^{2p-1}_j$.  Then $B$ preserves form parity  
$B_j = B_j^{ev} \oplus B_j^{odd}$ relative to the de Rham algebra written as a direct sum of even and odd forms. The self-adjoint first-order elliptic operators $B_j^{ev}$ and $B_j^{odd}$ are spectrally identical, one has
\begin{equation}\label{hjetaj} 
h_j := \Tr (\Pi_0[B_j^{ev}]) = \Tr (\Pi_0[B_j^{odd}])  = \frac{1}{2}\,\Tr (\Pi_0[B_j])
\end{equation} 
and 
$
\eta_j := \eta(B_j^{ev}, 0) = \eta(B_j^{odd}, 0) = \frac{1}{2}\,\eta(B_j, 0),
$
where $\Pi_0[S] \in \Fs_{\mbox{{\tiny $-\oo$}}}(\pd W_j)$ is the smoothing projection onto $\ker (S)$, and $\eta(S, 0)$ the $\eta$-invariant of an elliptic self-adjoint $\pdo$ $S$. Let 
\begin{equation}\label{pi0} 
\Pi_0^{ev} = \bigoplus_j \Pi_0[B_j^{ev}]  \hf \in  \Fs_{\mbox{{\tiny $-\oo$}}}(\pd W),
\end{equation}  
and likewise for $\Pi_0^{odd}$, and set
$
h := \Tr_{\mbox{{\tiny$\pd W$}}} (\Pi_0^{ev}) = \sum_j h_j $,  $\eta := \eta(B^{ev}, 0) = \sum_j\eta_j.
$
The APS projection is the order zero $\pdo$ projector
\begin{equation}\label{APS sum} 
 \Pi_\geq^{\mbox{{\tiny $ \partial W$}}} = \bigoplus_{j=1}^r \Pi_\geq^{\mbox{{\tiny $ \partial W_j$}}}  \ \in   \Fs_\Z(\pd W) := \bigoplus_{j=1}^r  \Psi^\Z( \pd W_j, \wedge T^*\pd W_j)
\end{equation}
where $\Pi_\geq^{\mbox{{\tiny $ \partial W_j$}}}$ is the orthogonal projection onto the span of eigenforms of $B_j$ with eigenvalue $\la\geq 0$. 
The Calder\'{o}n projection, on the other hand, 
$ C[\eth^{\mbox{{\tiny$W$}}}] \in \Fs_\Z(\pd W)$ is a projector onto the subspace $K(\eth^{\mbox{{\tiny$W$}}}) \< \O(\pd W)$ of  boundary sections which are restrictions to the boundary of interior solutions $\Ker(\eth^{\mbox{{\tiny$W$}}}) \< \O(W)$; the Poisson operator $\Kk[\eth^{\mbox{{\tiny$W$}}}]: \O(\pd W) \to \O(W)$  associated to $\eth^{\mbox{{\tiny$W$}}}$ restricts  in each Sobolev completion to a canonical isomorphism 
\begin{equation}\label{Poisson} 
K(\eth^{\mbox{{\tiny$W$}}}) \stackrel{\cong}{\to} \Ker(\eth^{\mbox{{\tiny$W$}}}) \hf \aand \ {\rm then} \hf 
C[\eth^{\mbox{{\tiny$W$}}}]  :=  \varrho\Kk[\eth^{\mbox{{\tiny$W$}}}], 
\end{equation}
where $ \varrho :  \O(W) \to \O(\pd W)$ is restriction to the boundary.  See for instance \S 7 of \cite{Gr}. 

Relative to an identification with its connected components $\pd W  = \pd W_1\sqcup\cdots\sqcup \pd W_n$ the projections may be written  as $n\x n$ block matrices: $\Pi_\geq^{\mbox{{\tiny $ \partial W$}}}$  is a diagonal direct sum of order zero $\pdos$ whilst the Calder\'{o}n projector $C[\eth^{\mbox{{\tiny$W$}}}]$ has order zero $\pdos$ along the diagonal and has non-zero off-diagonal smoothing operator terms.  
The crucial analytic fact is:
\begin{lem}\label{C-Pi} 
\begin{equation}\label{APS - C}  
C[\eth^{\mbox{{\tiny$W$}}}]  -  \Pi_\geq^{\mbox{{\tiny $ \partial W$}}} \ \in \Fs_{\mbox{{\tiny $-\oo$}}}(\pd W).
\end{equation}
\end{lem}
Proof:\; Since $\eth^{\mbox{{\tiny$W$}}}$ has the form $\s(du)(\pd_u + B_j)$ in a collar neighbourhood $U_i$  of each connected component $\pd W_i$, the argument in \cite{Sc cmp} (Prop. 2.2), or the more general argument of \cite{Gr} (Prop. 4.1), for the case for a single boundary readily adapts to the present multi-boundary context. 
\hfill
$\Box$\\[2mm]

The projection operators above are sensitive to orientation. For an oriented manifold $N$, let $N^-$ denote the manifold with orientation reversed. 
\begin{lem}\label{C-Pi orientation} 
$\Pi_\geq^{\mbox{{\tiny $\partial W^-$}}}  =  \Pi_\leq^{\mbox{{\tiny $ \partial W$}}}$
 is the projection onto the span of eigenforms with eigenvalue $\la \leq0$.  Likewise, $C[\eth^{\mbox{{\tiny$W$}}}]$ and $C[\eth^{\mbox{{\tiny$W^-$}}}]$ are complementary projections modulo smoothing operators.  
\end{lem}
Proof:\;  Reversing the orientation on $\pd W$ reverses the sign of the Riemannian volume form, and so the Hodge star $* \mto - *$. Thus $B^{2p}_j := (-1)^{k+p+1}(*d_j - d_j*)$  and $B^{2p-1}_j := (-1)^{k+p} (*d_j + d_j*)$ change sign, swapping negative and positive eigenvalues, which is the first assertion.  Since $\pd(W^-) = (\pd W)^-$, the statement for the Calder\'{o}n projection then follows from \eqref{APS - C}.
\hfill
$\Box$\\[2mm]

A representative $W$ for a bordism in $\mor_{\mbox{{\tiny $ {\rm Bord}_{4k}$}}}(M_0, M_1)$ comes with an orientation preserving diffeomorphism 
$
\kappa:  \pd W \to M_0^- \sqcup M_1.$
One has that   
$\kappa_\sharp(\Pi_\geq^{\mbox{{\tiny $ \partial W$}}} ),  \ \kappa_\sharp(C[\eth^{\mbox{{\tiny$W$}}}]) \in \Fs_\Z(M_0 \sqcup M_1)$
are order zero $\pdo$ projections, while
\begin{equation}\label{APS - kappa 2} 
\kappa_\sharp(C[\eth^{\mbox{{\tiny$W$}}}])  -  \kappa_\sharp(\Pi_\geq^{\mbox{{\tiny $ \partial W$}}} ) =  \kappa_\sharp(C[\eth^{\mbox{{\tiny$W$}}}]  -  \Pi_\geq^{\mbox{{\tiny $ \partial W$}}} ) \ \in \Fs_{\mbox{{\tiny $-\oo$}}}(M_0 \sqcup M_1)
\end{equation}
are smoothing operators. Also 
$\kappa_\sharp(\Pi_0^{ev})\ \in \Fs_{\mbox{{\tiny $-\oo$}}}(M_0 \sqcup M_1).$
To define a logarithm
$$\log^\sgn : \Nn \cob_{4k}\to \Fs_{\mbox{{\tiny $-\oo$}}}(\cob_{4k}^*)/[\Fs_{\mbox{{\tiny $-\oo$}}}(\cob_{4k}^*), \Fs_{\mbox{{\tiny $-\oo$}}}(\cob_{4k}^*)] $$ 
it is enough to specify it on 1-simplices 
\begin{equation*}
\log^\sgn_{\mbox{{\tiny $M_0 \sqcup M_1$}}} \ : \ \mor_{\mbox{{\tiny $ {\rm Bord}_{4k}$}}}(M_0, M_1) \to \Fs_{\mbox{{\tiny $-\oo$}}}(M_0 \sqcup M_1)/[\Fs_{\mbox{{\tiny $-\oo$}}}(M_0 \sqcup M_1), \Fs_{\mbox{{\tiny $-\oo$}}}(M_0 \sqcup M_1)].
\end{equation*}
Define
\begin{equation}\label{sgn log} 
\log^\sgn_{\mbox{{\tiny $M_0 \sqcup M_1$}}}  (\ol W) :=  \pi_{\mbox{{\tiny $M_0 \sqcup M_1$}}} \circ\kappa_\sharp\l(C[\eth^{\mbox{{\tiny$W$}}}] - \Pi_\geq^{\mbox{{\tiny $ \partial W$}}}  + \Pi_0^{ev} \r)
\end{equation}

\hf\hf --- equal to the sum of order zero $\pdo$ projections in $\Fs^{\,0}_{\mbox{{\tiny $\Z, -\oo$}}}(M_0 \sqcup M_1)$ ---
\begin{equation*}
=  \pi_{\mbox{{\tiny $M_0 \sqcup M_1$}}} \circ \kappa_\sharp (C[\eth^{\mbox{{\tiny$W$}}}]) - \pi_{\mbox{{\tiny $M_0 \sqcup M_1$}}} \circ \kappa_\sharp(\Pi_\geq^{\mbox{{\tiny $ \partial W$}}} ) + \pi_{\mbox{{\tiny $M_0 \sqcup M_1$}}} \circ \kappa_\sharp(\Pi_0^{ev}).
\end{equation*}
From \eqref{trace diagram theta mn} and \eqref{phi sharp 2} 
\begin{equation}\label{sgn log 2} 
\log^\sgn_{\mbox{{\tiny $M_0 \sqcup M_1$}}}  (\ol W) =  \vartheta_{\mbox{{\tiny $\partial W, M_0 \sqcup M_1$}}} \circ \pi_{\mbox{{\tiny $ \partial W$}}} \l(C[\eth^{\mbox{{\tiny$W$}}}] - \Pi_\geq^{\mbox{{\tiny $ \partial W$}}}  + \Pi_0^{ev}\r).
\end{equation}

\vskip 3mm

\begin{prop}\label{sign char}  The right-hand side of 
\eqref{sgn log} depends only on the (oriented) bordism class  $\ol W$ of $W$ (independent of $g_W$)  and has log-character 
\begin{equation}\label{sgn log char} 
\wt\Tr_{\mbox{{\tiny $M_0 \sqcup M_1$}}}(\log^\sgn_{\mbox{{\tiny $M_0 \sqcup M_1$}}}  \ol W) =  \sgn(W).
\end{equation}
\end{prop} 

For use here and elsewhere, we note the following lemma:
\begin{lem}\label{Hilbert trace index} 
Let $H = H_+ \oplus H_-$ be a Hilbert space polarised by infinite-dimensional subspaces $H_\pm$, and let $\Pi_\pm$ be the  orthogonal projections with ranges $H_\pm$. Let $P_0, P_1$ be projections on $H$ with  $P_j - \Pi_+$ of trace-class ($j=0,1$) on $H$.  Let $W_j:=\ran(P_j)\<H$, and 
let $\ind_{\mbox{{\tiny $W_0, W_1$}}}   a$ denote the index of a Fredholm operator $a: W_0 \to W_1$. 
Then $P_0 - P_1$ is trace class on $H$ and $P_1P_0 : W_0 \to W_1$ is a Fredholm operator, and one has  
\begin{equation}\label{ind=tr} 
\ind_{\mbox{{\tiny $W_0, W_1$}}} (P_1P_0) = \Tr_H (P_0 - P_1).
\end{equation}
\end{lem}
Proof:\;  Follows in a straightforward way using the methods of \S 7.1 of \cite{PrSe}.
\hfill
$\Box$\\[2mm]

Proof of \propref{sign char}:\; Let $\eth^{\mbox{{\tiny$W$}}}_{\mbox{{\tiny$\geq$}}}$ be the APS boundary value problem \cite{APS1}. Thus, $\eth^{\mbox{{\tiny$W$}}}_{\mbox{{\tiny$\geq$}}} = \eth^{\mbox{{\tiny$W$}}}$ with domain restricted to those sections $s\in \O^+(W)$ with  $\Pi_{\mbox{{\tiny$\geq$}}}^{\mbox{{\tiny $ \partial W$}}} (s_{|{\mbox{{\tiny $ \partial W$}}}}) =0$. Then, in the notation of \lemref{Hilbert trace index}, 
\begin{equation}\label{ind PK} 
\ind \eth^{\mbox{{\tiny$W$}}}_{\mbox{{\tiny$\geq$}}} = 
\ind_{\mbox{{\tiny $K(\eth^{\mbox{{\tiny$W$}}}_{\mbox{{\tiny$\geq$}}}), \ran(\Pi_\geq^{\mbox{{\tiny $ \partial W$}}} )$}}}  \l(\Pi_{\mbox{{\tiny$\geq$}}}^{\mbox{{\tiny $ \partial W$}}} \circ C(\eth^{\mbox{{\tiny$W$}}}_{\mbox{{\tiny$\geq$}}} ) \r)
\end{equation}
with $ K(\eth^{\mbox{{\tiny$W$}}}_{\mbox{{\tiny$\geq$}}})$ in \eqref{Poisson} viewed as a closed subspace of the Hilbert space $H^{\mbox{{\tiny $ \partial W$}}}$ of $L^2$ boundary sections polarised with $H^{\mbox{{\tiny $ \partial W$}}}_+ =  \ran(\Pi_\geq^{\mbox{{\tiny $ \partial W$}}} )$, $H^{\mbox{{\tiny $ \partial W$}}}_- =  \ran(\Pi_<^{\mbox{{\tiny $ \partial W$}}})$ (the identity \eqref{ind PK}  for Dirac-type  operators is well known, see for instance \cite{BoWo}, \cite{Sc cmp}). With $h$ and $\eta$ defined following \eqref{pi0} and $L(w)$ the Hirzebruch $L$-polynomial in the Pontryagin forms, the APS signature theorem gives the first two equalities in
\begin{eqnarray*}
\sgn(W) & \stackrel{\mbox{{\tiny\cite{APS1}, Thm 4.14}}}{=} & \int_W L(w) - \eta \\
& \stackrel{\mbox{{\tiny\cite{APS1}, eqn 4.7}}}{=} &\ind (\eth^{\mbox{{\tiny$W$}}}_{\mbox{{\tiny$\geq$}}}) + h  \\
& \stackrel{\mbox{{\tiny\eqref{ind PK}}}}{=} &\ind_{\mbox{{\tiny $K(\eth^{\mbox{{\tiny$W$}}}_{\mbox{{\tiny$\geq$}}}), \ran(\Pi_\geq^{\mbox{{\tiny $ \partial W$}}} )$}}}  \l(\Pi_{\mbox{{\tiny$\geq$}}}^{\mbox{{\tiny $ \partial W$}}} \circ C[\eth^{\mbox{{\tiny$W$}}}_{\mbox{{\tiny$\geq$}}}]\r)
+ \Tr_{\mbox{{\tiny$\pd W$}}}(\Pi_0^{ev})  \\ 
& \stackrel{\mbox{{\tiny\eqref{ind=tr}}}}{=} & \Tr_{\mbox{{\tiny $\pd W$}}}(C[\eth^{\mbox{{\tiny$W$}}}] - \Pi_{\mbox{{\tiny$\geq$}}}^{\mbox{{\tiny $ \partial W$}}}) + \Tr_{\mbox{{\tiny$\pd W$}}}(\Pi_0^{ev})   \\[2mm]
& = & \Tr_{\mbox{{\tiny $\pd W$}}}(C[\eth^{\mbox{{\tiny$W$}}}] - \Pi_{\mbox{{\tiny$\geq$}}}^{\mbox{{\tiny $ \partial W$}}} + \Pi_0^{ev})   \\
&= & \Tr_{\mbox{{\tiny $M_0 \sqcup M_1$}}}(\kappa_\sharp(C[\eth^{\mbox{{\tiny$W$}}}] - \Pi_{\mbox{{\tiny$\geq$}}}^{\mbox{{\tiny $ \partial W$}}} + \Pi_0^{ev}))   \\
& \stackrel{\mbox{{\tiny\eqref{commuting trace reduced}}}}{=} &  \wt\Tr_{\mbox{{\tiny $M_0 \sqcup M_1$}}}(\log^\sgn_{\mbox{{\tiny $M_0 \sqcup M_1$}}} \ol W).
\end{eqnarray*}

The character $\wt\Tr_{\mbox{{\tiny $M_0 \sqcup M_1$}}}(\log^\sgn_{\mbox{{\tiny $M_0 \sqcup M_1$}}} \ol W)\in\Z$ is thus an oriented-homotopy  invariant of $W$. Since $\wt\Tr_{\mbox{{\tiny $M_0 \sqcup M_1$}}} : \Fs_{\mbox{{\tiny $-\oo$}}}(M_0 \sqcup M_1)/[\Fs_{\mbox{{\tiny $-\oo$}}}(M_0 \sqcup M_1), \Fs_{\mbox{{\tiny $-\oo$}}}(M_0 \sqcup M_1)] \stackrel{\cong}{\to} \C$ is  a linear isomorphism  by \lemref{Tr unique}, $\log^\sgn _{\mbox{{\tiny $M_0 \sqcup M_1$}}}  \ol W$ is hence a homotopy invariant of the manifold $W$; that is, with $\simeq_O$ indicating oriented homotopy equivalence,  $$W\simeq_O W\pr \To \sgn W = \sgn W\pr \To \wt\Tr_{\mbox{{\tiny $M_0 \sqcup M_1$}}}(\log^\sgn_{\mbox{{\tiny $M_0 \sqcup M_1$}}} \ol W  - \log^\sgn_{\mbox{{\tiny $M_0 \sqcup M_1$}}} \ol W\pr) =0 $$ $$ \To  \log^\sgn_{\mbox{{\tiny $M_0 \sqcup M_1$}}} \ol W  = \log^\sgn_{\mbox{{\tiny $M_0 \sqcup M_1$}}} \ol W\pr \hf \iin \ \Fs_{\mbox{{\tiny $-\oo$}}}(M_0 \sqcup M_1)/[\Fs_{\mbox{{\tiny $-\oo$}}}(M_0 \sqcup M_1), \Fs_{\mbox{{\tiny $-\oo$}}}(M_0 \sqcup M_1)].$$
In particular, the logarithm is an invariant of the bordism class of $W$ in $\mor_{\mbox{{\tiny $ {\rm Bord}_{4k}$}}}(M_0, M_1)$, and independent of any choice of Riemannian metric on $W$.
\begin{flushright}
$\Box$
\end{flushright}

It is useful to note: 
\begin{lem}\label{B even odd}
$\log^\sgn_{\mbox{{\tiny $M_0 \sqcup M_1$}}}  (\ol W)$ in \eqref{sgn log}, or \eqref{sgn log 2}, is unchanged if $B^{ev}$ is replaced by $B^{odd}$
\end{lem}
Proof:\;  The difference is $\pi_{\mbox{{\tiny $M_0 \sqcup M_1$}}} \circ\kappa_\sharp\l( \Pi_0^{ev} - \Pi_0^{odd} \r)$ which has character 
$$\wt\Tr_{\mbox{{\tiny $M_0 \sqcup M_1$}}}(\pi_{\mbox{{\tiny $M_0 \sqcup M_1$}}} \circ\kappa_\sharp\l( \Pi_0^{ev} - \Pi_0^{odd}\r)) = \Tr_{\mbox{{\tiny $M_0 \sqcup M_1$}}}( \Pi_0^{ev} - \Pi_0^{odd})$$ which by \eqref{hjetaj} is zero. Since  $\wt\Tr_{\mbox{{\tiny $M_0 \sqcup M_1$}}}$ is a  isomorphism, the assertion follows.
$\hfill\Box$

We may therefore better write
\begin{eqnarray*}
\log^\sgn_{\mbox{{\tiny $M_0 \sqcup M_1$}}}  (\ol W) & = & \pi_{\mbox{{\tiny $M_0 \sqcup M_1$}}} \circ\kappa_\sharp\l(C[\eth^{\mbox{{\tiny$W$}}}] - \Pi_\geq^{\mbox{{\tiny $ \partial W$}}}  + U^{\mbox{{\tiny $ \partial W$}}}\r) \label{sgn log 3} \\[2mm]
& = & \vartheta_{\mbox{{\tiny $\partial W, M_0 \sqcup M_1$}}} \circ \pi_{\mbox{{\tiny $ \partial W$}}} \l(C[\eth^{\mbox{{\tiny$W$}}}] - \Pi_\geq^{\mbox{{\tiny $ \partial W$}}}  + U^{\mbox{{\tiny $ \partial W$}}}\r)\label{sgn log 4} 
\end{eqnarray*}
with $U^{\mbox{{\tiny $ \partial W$}}}$ denoting either of the projections; this flexibility is important later.

The principal task at hand is to show log-additivity: 

\begin{thm}\label{log sgn additivity} With respect to composition of bordisms 
$$  \mor_{\mbox{{\tiny $ {\rm Bord}_{4k}$}}}(M_0, M_1)  \x  \mor_{\mbox{{\tiny $ {\rm Bord}_{4k}$}}}(M_1, M_2)  \to  \mor_{\mbox{{\tiny $ {\rm Bord}_{4k}$}}}(M_0, M_2), \ \  (\ol W_0, \ol W_1) \mto \ol W_0\cup\ol W_1,$$
one has in $\Fs_{\mbox{{\tiny $-\oo$}}}(M_0 \sqcup M_1 \sqcup M_2) / [\Fs_{\mbox{{\tiny $-\oo$}}}(M_0 \sqcup M_1 \sqcup M_2),  \Fs_{\mbox{{\tiny $-\oo$}}}(M_0 \sqcup M_1 \sqcup M_2)]$
\begin{equation}\label{log sgn add} 
\wt\eta_{\mbox{{\tiny $M_1$}}}\log^\sgn_{\mbox{{\tiny $M_0\sqcup M_2$}}}(\ol W_0\cup\ol W_1) 
\ =\ \wt\eta_{\mbox{{\tiny $M_2$}}}\log^\sgn_{\mbox{{\tiny $M_0\sqcup M_1$}}}(\ol W_0) 
\ + \ \wt\eta_{\mbox{{\tiny $M_0$}}}\log^\sgn_{\mbox{{\tiny $M_1\sqcup M_2$}}}(\ol W_1).
\end{equation}
\end{thm}

\vskip 2mm 
Applying the trace $\wt\Tr_{\mbox{{\tiny $M_0 \sqcup M_1 \sqcup M_2$}}}$ to \eqref{log sgn add},  one has from \eqref{sgn log char}:

\begin{cor}\label{only on W}
\begin{equation}\label{log signature}
\sgn(W\cup_{M_1} W\pr) =  \sgn(W) + \sgn(W\pr).
\end{equation}
\end{cor}

 \eqref{log signature} was originally observed by Novikov (c1967) 
\footnote{Contrasting with (Wall) non-additivity of the signature for higher codimension partitions \cite{Wall}.} and proved 
 for closed $W\cup_{M_1} W\pr$ in \cite{ASIII}. 
 
\begin{cor}\label{only on W} 
$\log^\sgn_{\mbox{{\tiny $M_0\sqcup M_1$}}}(\ol W_0)$ is independent of the boundary diffeomorphism $\kappa$, and so depends only on the oriented diffeomorphism class of $W$  (in fact, homotopy class). $\log^\sgn_{\mbox{{\tiny $M_0\sqcup M_2$}}}(\ol W_0\cup \ol W_1)$ is independent of the gluing diffeomorphism $\phi$ between the identified  boundary components  of $W_0\in \ol W_0$ and $W_1\in \ol W_1$ used to form $\ol W_0\cup \ol W_1 := \ol{W_0\cup_\phi W_1}$.  The same statements hold for  $\sgn(W_0)$ and $\sgn(W_0 \cup_\phi W_1)$. 
\end{cor}

The proof of \thmref{log sgn additivity} will occupy the remainder of this section. 

\begin{prop}\label{log sgn additivity 2} The equality \eqref{log sgn add}  holds if 
\begin{equation}\label{log sgn add 2} 
\wt\eta_{\mbox{{\tiny $M_1 \sqcup M_1$}}}\log^\sgn_{\mbox{{\tiny $M_0\sqcup M_2$}}}(\ol W_0\cup \ol W_1) 
\ =\ \wt\eta_{\mbox{{\tiny $M_1 \sqcup M_2$}}}\log^\sgn_{\mbox{{\tiny $M_0\sqcup M_1$}}}(\ol W_0) 
\ + \ \wt\eta_{\mbox{{\tiny $M_0\sqcup M_1$}}}\log^\sgn_{\mbox{{\tiny $M_1\sqcup M_2$}}}(\ol W_1) 
\end{equation}
holds in $\Fs_{\mbox{{\tiny $-\oo$}}}(M_0 \sqcup M_1 \sqcup M_1 \sqcup M_2)/[\Fs_{\mbox{{\tiny $-\oo$}}}(M_0 \sqcup M_1 \sqcup M_1 \sqcup M_2),  \Fs_{\mbox{{\tiny $-\oo$}}}(M_0 \sqcup M_1 \sqcup M_1 \sqcup M_2)].$
\end{prop}

Proof:\;   
$$
\wt\Tr_{\mbox{{\tiny $M_0 \sqcup M_1 \sqcup M_1 \sqcup M_2$}}}(\wt\eta_{\mbox{{\tiny $M_1 \sqcup M_1$}}}\log^\sgn_{\mbox{{\tiny $M_0\sqcup M_2$}}}(\ol W_0\cup\ol W_1))   \  \stackrel{\eqref{tilde tracial Finfty}}{=} \  \wt\Tr_{\mbox{{\tiny $M_0  \sqcup M_2$}}}(\log^\sgn_{\mbox{{\tiny $M_0\sqcup M_2$}}}(\ol W_0\cup \ol W_1))$$
$$ \stackrel{\eqref{tilde tracial Finfty}}{=}    \ \   \wt\Tr_{\mbox{{\tiny $M_0  \sqcup M_1 \sqcup M_2$}}}(\wt\eta_{\mbox{{\tiny $M_1$}}}\log^\sgn_{\mbox{{\tiny $M_0\sqcup M_2$}}}(\ol W_0\cup\ol W_1)),$$
and, similarly, 
$$\wt\Tr_{\mbox{{\tiny $M_0 \sqcup M_1 \sqcup M_1 \sqcup M_2$}}}( \wt\eta_{\mbox{{\tiny $M_1 \sqcup M_2$}}}\log^\sgn_{\mbox{{\tiny $M_0\sqcup M_1$}}}(\ol W_0)) =  \wt\Tr_{\mbox{{\tiny $M_0 \sqcup M_1 \sqcup M_2$}}}( \wt\eta_{\mbox{{\tiny $M_2$}}}\log^\sgn_{\mbox{{\tiny $M_0\sqcup M_1$}}}(\ol W_0)),$$
 $$\wt\Tr_{\mbox{{\tiny $M_0 \sqcup M_1 \sqcup M_1 \sqcup M_2$}}}( \wt\eta_{\mbox{{\tiny $M_0\sqcup M_1$}}}\log^\sgn_{\mbox{{\tiny $M_1\sqcup M_2$}}}(\ol W_1)) =   \wt\Tr_{\mbox{{\tiny $M_0 \sqcup M_1 \sqcup M_2$}}}( \wt\eta_{\mbox{{\tiny $M_0$}}}\log^\sgn_{\mbox{{\tiny $M_1\sqcup M_2$}}}(\ol W_1)).$$
Hence,  if \eqref{log sgn add 2} holds, $\wt\Tr_{\mbox{{\tiny $M_0 \sqcup M_1 \sqcup M_2$}}}$ evaluated on 
$$\wt\eta_{\mbox{{\tiny $M_1$}}}\log^\sgn_{\mbox{{\tiny $M_0\sqcup M_2$}}}(\ol W_0\cup\ol W_1) - \wt\eta_{\mbox{{\tiny $M_2$}}}\log^\sgn_{\mbox{{\tiny $M_0\sqcup M_1$}}}(\ol W_0) -  \wt\eta_{\mbox{{\tiny $M_0$}}}\log^\sgn_{\mbox{{\tiny $M_1\sqcup M_2$}}}(\ol W_1)$$
is zero. Since $\wt\Tr_{\mbox{{\tiny $M_0 \sqcup M_1 \sqcup M_2$}}}$  is  from \eqref{trace reduced} a linear isomorphism, \eqref{log sgn add}  follows.
$\hfill\Box$

\vskip 5mm

\corref{only on W} allows one to  work with the geometric boundary  of a representative $W_0$ of $\ol W \in \mor_{\mbox{{\tiny $ {\rm Bord}_{4k}$}}}(M_0, M_1)$, rather than $M_0, M_1$. Thus, $\pd W_0 = X^-_0 \sqcup X_1$ along with orientation preserving diffeomorphisms $\a_{\mbox{{\tiny$\pd W_0$}}} : X_0 \to M_0$ and $\b_{\mbox{{\tiny$\pd W_0$}}} : X_1 \to M_1$. Likewise,  $W_1\in \ol W_1\in \mor_{\mbox{{\tiny $ {\rm Bord}_{4k}$}}}(M_1, M_2) $ has $\pd W_1 =   Y^-_1 \sqcup Y_2$ and oriented diffeomorphisms $\a_{\mbox{{\tiny$\pd W_1$}}} : Y_1 \to M_1$ and $\b_{\mbox{{\tiny$\pd W_1$}}} : Y_2 \to M_2$. Let  $\phi =  \a_{\mbox{{\tiny$\pd W_1$}}} \ii \circ \b_{\mbox{{\tiny$\pd W_0$}}}  : X_1 \stackrel{\cong}{\to} Y_1 $. The space $W_0\cup_\phi W_1$  formed from $W_0$ and $W_1$ by identifying $x\in X_1$ with $\phi(x)\in Y_1$ has a smooth manifold structure compatible with those of $W_0$ and $W_1$ which is unique modulo oriented diffeomorphisms which fix $M_0, \phi(X_1) = Y_1$ and $M_2$. Then $\ol W_0\cup\ol  W_1 : = \ol{W_0\cup_\phi W_1} \in \mor_{\mbox{{\tiny $ {\rm Bord}_{4k}$}}}(M_0, M_2) $ is the equivalence class of $W_0\cup_\phi W_1$ modulo such diffeomorphisms compatible with $\a_{\mbox{{\tiny$\pd W_0$}}}$ and $\b_{\mbox{{\tiny$\pd W_1$}}}$. 
One has, further, the closed oriented hypersurface 
$
N =\{[x] \ | \ x\in X_1\} \ \<  \  W_0\cup_\phi W_1$
with $[x]$ the equivalence class in the identification space  $W_0\cup_\phi W_1$. We may choose a  choose a Riemannian metric on $W_0\cup_\phi W_1$ which is isometric to a product in some collar neighbourhood $U \cong (-1,1) \x N$ of $N$ in $W_0\cup_\phi W_1$, with $N$ identified with $\{0\}\x N \< U$. 
Define, then, 
\begin{equation*}\label{log W0}
\log_{\mbox{{\tiny $X_0 \sqcup X_1$}}}  (\ol W_0)  := \pi_{\mbox{{\tiny $X_0 \sqcup X_1$}}} \l(C[\eth^{\mbox{{\tiny$W_0$}}}] - \Pi_\geq^{\mbox{{\tiny $ X^-_0 \sqcup X_1$}}}  + U^{\mbox{{\tiny $ X_0 \sqcup X_1$}}}\r), 
\end{equation*} 

\begin{equation*}\label{log W1}
\log_{\mbox{{\tiny $Y_1 \sqcup Y_2$}}}  (\ol W_1)  := \pi_{\mbox{{\tiny $Y_1 \sqcup Y_2$}}} \l(C[\eth^{\mbox{{\tiny$W_1$}}}] - \Pi_\geq^{\mbox{{\tiny $ Y^-_1 \sqcup Y_2$}}}  + U^{\mbox{{\tiny $ Y_1 \sqcup Y_2$}}}\r), 
\end{equation*} 

\begin{equation*}\label{log W0W1}
\log_{\mbox{{\tiny $X_0 \sqcup Y_2$}}}  (\ol W_0\cup\ol W_1)  := \pi_{\mbox{{\tiny $X_0 \sqcup Y_2$}}} \l(C[\eth^{\mbox{{\tiny$W_0\cup_\phi W_1$}}}] - \Pi_\geq^{\mbox{{\tiny $ X^-_0 \sqcup Y_2$}}}  + U^{\mbox{{\tiny $ X_0 \sqcup Y_2$}}}\r). \\[2mm]
\end{equation*} 
In terms other than $\Pi_\geq$ the orientation is not felt and so is not indicated. 

\begin{prop}\label{log sgn additivity 3} The equality  \eqref{log sgn add 2}  holds if 
\begin{equation}\label{log sgn add 3} 
\wt\eta_{\mbox{{\tiny $X_1 \sqcup Y_1$}}}\log_{\mbox{{\tiny $X_0 \sqcup Y_2$}}}  (\ol W_0\cup\ol W_1) 
\ =\ \wt\eta_{\mbox{{\tiny $Y_1 \sqcup Y_2$}}}\log_{\mbox{{\tiny $X_0 \sqcup X_1$}}}  (\ol W_0)    
\ + \ \wt\eta_{\mbox{{\tiny $X_0\sqcup X_1$}}}\log_{\mbox{{\tiny $Y_1 \sqcup Y_2$}}}  (\ol W_1) 
\end{equation}
holds in $\Fs_{\mbox{{\tiny $-\oo$}}}(X_0 \sqcup X_1 \sqcup Y_1 \sqcup Y_2)/[\Fs_{\mbox{{\tiny $-\oo$}}}(X_0 \sqcup X_1 \sqcup Y_1 \sqcup Y_2),  \Fs_{\mbox{{\tiny $-\oo$}}}(X_0 \sqcup X_1 \sqcup Y_1 \sqcup Y_2)].$ 
\end{prop}

Proof:\; Let $V_j, Z_j, M, N \in \ob(\cob_n)$ with $V_j$ and $Z_j$ diffeomorphic and $M$ and $N$ diffeomorphic. Let $V := V_1 \sqcup \cdots \sqcup V_m$ and $Z := Z_1 \sqcup \cdots \sqcup Z_m$.
 By \eqref{phi sharp 2}, there are then canonical identifications  $\theta_{V,Z}:  \Fs_\Pi(V)\to  \Fs_\Pi(Z)$ and $\vartheta_{\mbox{{\tiny $V_{\mbox{{\tiny $N$}}}, Z_{\mbox{{\tiny $M$}}}$}}}:  \Fs_\Pi(V_{\mbox{{\tiny $N$}}})  \to  \Fs_\Pi(Z_{\mbox{{\tiny $M$}}})$, where
$$V_{\mbox{{\tiny $N$}}} := V_1 \sqcup \cdots X_{k-1} \sqcup N \sqcup X_k \sqcup \cdots \sqcup V_m, \hf Z_M := Z_1 \sqcup \cdots Z_{k-1} \sqcup M \sqcup Z_k \sqcup \cdots \sqcup Z_m.$$
Moreover, the following diagram is easily seen to commute
\begin{equation}\label{diagram theta eta} 
\begin{array}{rcl}
  \Fs_\Pi(V_{\mbox{{\tiny $N$}}}) \ &   \stackrel{\vartheta_{\mbox{{\tiny $V_{\mbox{{\tiny $N$}}}, Z_{\mbox{{\tiny $M$}}}$}}}}{\longrightarrow}   & \  \Fs_\Pi(Z_{\mbox{{\tiny $M$}}})   \\
    &  &   \\
  \uparrow\wt\eta^k_{\mbox{{\tiny $N$}}}  &    & \hf\uparrow\wt\eta^k_{\mbox{{\tiny $M$}}} \\
    &  &   \\
 \Fs_\Pi(V)    & \stackrel{\vartheta_{\mbox{{\tiny $V, Z$}}}}{\longrightarrow}   &   \Fs_\Pi(Z)
\end{array}.
\end{equation}  
Hence, taking $M:= M_1 \sqcup M_2, N := Y_1 \sqcup Y_2,  V:= X_0 \sqcup X_1,  Z:= M_0 \sqcup M_1,V_M : = X_0 \sqcup X_1\sqcup Y_1 \sqcup Y_2, Z_M:= M_0 \sqcup M_1\sqcup M_1 \sqcup M_2,$
one has
\begin{eqnarray*}
\wt\eta_{\mbox{{\tiny $M_1 \sqcup M_2$}}}\log^\sgn_{\mbox{{\tiny $M_0\sqcup M_1$}}}(\ol W_0)  & = & 
\wt\eta_{\mbox{{\tiny $M_1 \sqcup M_2$}}}\circ\vartheta_{\mbox{{\tiny $X_0 \sqcup X_1, M_0 \sqcup M_1$}}} \log_{\mbox{{\tiny $X_0 \sqcup X_1$}}}  (\ol W_0) \\ 
& =& 
\vartheta_{\mbox{{\tiny $X_0 \sqcup X_1\sqcup Y_1 \sqcup Y_2, M_0 \sqcup M_1\sqcup M_1 \sqcup M_2$}}} \l( \wt\eta_{\mbox{{\tiny $Y_1 \sqcup Y_2$}}}\log_{\mbox{{\tiny $X_0 \sqcup X_1$}}}  (\ol W_0)\r),
\end{eqnarray*}
$$\wt\eta_{\mbox{{\tiny $M_0 \sqcup M_1$}}}\log^\sgn_{\mbox{{\tiny $M_1\sqcup M_2$}}}(\ol W_1) = 
\vartheta_{\mbox{{\tiny $X_0 \sqcup X_1\sqcup Y_1 \sqcup Y_2, M_0 \sqcup M_1\sqcup M_1 \sqcup M_2$}}} \l( \wt\eta_{\mbox{{\tiny $X_0 \sqcup X_1$}}}\log_{\mbox{{\tiny $Y_1 \sqcup Y_2$}}}  (\ol W_1)\r)$$
$$\wt\eta_{\mbox{{\tiny $M_1 \sqcup M_1$}}}\log^\sgn_{\mbox{{\tiny $M_0\sqcup M_2$}}}(\ol W_0 \cup \ol W_1) = 
\vartheta_{\mbox{{\tiny $X_0 \sqcup X_1\sqcup Y_1 \sqcup Y_2, M_0 \sqcup M_1\sqcup M_1 \sqcup M_2$}}} \l( \wt\eta_{\mbox{{\tiny $X_1 \sqcup Y_1$}}}\log_{\mbox{{\tiny $X_0 \sqcup Y_2$}}}  (\ol W_0 \cup \ol W_1)\r).$$
Hence
$ \eqref{log sgn add 2} \ = \   \underbrace{\vartheta_{\mbox{{\tiny $X_0 \sqcup X_1\sqcup Y_1 \sqcup Y_2, M_0 \sqcup M_1\sqcup M_1 \sqcup M_2$}}}}_{\mbox{{\tiny linear isomorphism}}} \l(\eqref{log sgn add 3}\r).$
$\hfill\Box$

\begin{prop}\label{log sgn additivity final} \ \\[2mm]
The equality  \eqref{log sgn add 3} holds.
\end{prop}

Proof:\; It is convenient to take $W\in \ol W_0$ and $W\pr\in \ol W_1$ by cutting $W_0 \cup_\phi W_1 \in \ol W_0 \cup \ol W_1$ along the hypersurface N: let
$W := (W_0 \cup_\phi W_1)\bsh (W_1\bsh N), W\pr := (W_0 \cup_\phi W_1)\bsh (W_0\bsh N).$
Set $X:= X_0, Y:= Y_2.$ Then  
\begin{equation}\label{idents} 
\pd W = X^- \sqcup N, \hf \pd W\pr = N^- \sqcup Y, \hf  X_1= N = Y_1.
\end{equation}
From the sequences of inclusions
$N \rightrightarrows W \sqcup W\pr  \to  W \cup_\phi W\pr$
one has the Mayer-Vietoris type sequence 
$
0\to \O^*(W_0 \cup_\phi W_1) \to \O^*(W) \oplus \O^*(W\pr) \to \O^*(N)$
in which the first map is signed restriction of a form 
$
\omega \mto (\omega_{|W}, -\,\omega_{|W\pr})$
and the second the sum of the boundary restrictions
$
(\sigma, \sigma\pr) \mto \sigma_{|N} + \sigma\pr_{|N}
$
(`restriction' means $\sigma_{|N} := j_N^*(\sigma)$ for $j_N:N\hookrightarrow W$ the inclusion, and so on). We assume for now that at least one of $W$ and $W\pr$ has disconnected boundary. Then the non-exact sequence Mayer Vietori sequence becomes exact on restriction to the kernels 
\begin{equation}\label{Mayer Vietoris 3} 
0\to \Ker(\eth^{\mbox{{\tiny $W \cup_\phi W\pr $}}}) \to  \Ker(\eth^{\mbox{{\tiny $W$}}}) \oplus \Ker(\eth^{\mbox{{\tiny $W\pr$}}}) \to \O^*(N) \to 0,
\end{equation}
by observing that $\Ker(\eth^{\mbox{{\tiny $W \cup_\phi W\pr $}}}) $ is the kernel of the map  $\Ker(\eth^{\mbox{{\tiny $W$}}}) \oplus \Ker(\eth^{\mbox{{\tiny $W\pr$}}})\to \O^*(N)$. But in a open set  $U = (-1,1) \x Y$, with $Y$ an odd-dimensional compact boundaryless manifold, the Riemannian metric can be chosen to be a product metric $g_{|U}  = du^2 + g_{{\mbox{{\tiny$Y$}}}}$, and so that $\eth^{\mbox{{\tiny $U$}}} = (du\mbox{{\small$\wedge$}} + i_{du})\l(\pd_u + D_{{\mbox{{\tiny$Y$}}}}\r)$ relative to the (self-adjoint) signature operator $\eth^{{\mbox{{\tiny$Y$}}}}$ on $Y$. This implies any solution $\psi$ to $\eth^{\mbox{{\tiny $U$}}}$ has the form 
$\psi(u,y) = \sum_k e^{-\la_k u}\psi_k(0)\phi_k(y)$ for a spectral resolution $(\la_k, \phi_k)$ of $\eth^{{\mbox{{\tiny$Y$}}}}$. The metric on $W \cup_{\mbox{{\tiny $N$}}} W\pr$ may be chosen to be a product in a tubular neighbourhood $(-1,1) \x N$ of the partitioning hypersurface $N$. Hence, matching of higher normal derivatives along $N$ of elements of $\Ker(\eth^{\mbox{{\tiny $W$}}})$ and $\Ker(\eth^{\mbox{{\tiny $W\pr$}}})$ follows from their zeroeth order matching pointwise along $N$ (with a change of sign taking into account the sign of $u$ in $(-1,1)$).

In view of the isomorphism \eqref{Poisson}, restricting solutions to the boundaries of the manifolds $W$ and $W\pr$ refines \eqref{Mayer Vietoris 3} to an exact sequence of maps on boundary sections
\begin{equation}\label{Mayer Vietoris 4}  
0\to K(\eth^{\mbox{{\tiny $W \cup_\phi W\pr $}}}) \to  K(\eth^{\mbox{{\tiny $W$}}}) \oplus K(\eth^{\mbox{{\tiny $W\pr$}}}) \to \O^*(N) \to 0.
\end{equation}
Let  $H^N$ be the space of forms $\O(N)$, or in the following can be taken to be its $L^2$ completion, on $N$.
The sequence \eqref{Mayer Vietoris 4} fits into a diagram 
\begin{equation}\label{diagram}  
\begin{array}{ccccccccc}
0 &  \to & K(\eth^{\mbox{{\tiny $W \cup_\phi W\pr $}}}) & \to  & K(\eth^{\mbox{{\tiny $W$}}})\oplus K(\eth^{\mbox{{\tiny $W\pr$}}}) &  \to & H^N & \to & 0 \\[4mm]
 &   & \  \ \  \downarrow G_0 &   & \hf  \downarrow G_1  &  & \  \ \ \downarrow \mbox{{\tiny $  id$}} &  &  \\[5mm]
0 &  \to & \ran(\Pi_>^{\mbox{{\tiny $\pd(W \cup_\phi W\pr) $}}}\oplus U^{\mbox{{\tiny $\pd(W \cup_\phi W\pr) $}}})    & \to  &  \begin{array}{c} \ran(\Pi_>^{\mbox{{\tiny $\pd W$}}}\oplus U^{\mbox{{\tiny $\pd W$}}}) \\[1mm] \oplus  \\[1mm] \ran(\Pi_>^{\mbox{{\tiny $\pd W\pr$}}}\oplus U^{\mbox{{\tiny $\pd W\pr$}}}) \end{array} & \to &  H^N & \to & 0\\[2mm]
&&&&&&&
 \end{array}
\end{equation}
where  in $\Psi^0(X\sqcup N \sqcup Y)$
\begin{eqnarray} 
G_0 & = & (\Pi_>^{\mbox{{\tiny $\pd(W \cup_\phi W\pr) $}}}\oplus U^{\mbox{{\tiny $\pd(W \cup_\phi W\pr)$}}})\circ C[\eth^{\mbox{{\tiny$W_0 \cup_\phi W_1$}}}], \nonumber \\[2mm]
G_1 & = & (\Pi_>^{\mbox{{\tiny $\pd W$}}}\oplus U^{\mbox{{\tiny $\pd W$}}})\circ C[\eth^{\mbox{{\tiny$W$}}}] \ \oplus \ (\Pi_>^{\mbox{{\tiny $\pd W\pr$}}}\oplus U^{\mbox{{\tiny $\pd W\pr$}}})\circ C[\eth^{\mbox{{\tiny$W\pr$}}}],\label{G1a}  \\[2mm]
& = & \l((\Pi_>^{\mbox{{\tiny $\pd W$}}}\oplus U^{\mbox{{\tiny $\pd W$}}}) \oplus \ (\Pi_>^{\mbox{{\tiny $\pd W\pr$}}}\oplus U^{\mbox{{\tiny $\pd W\pr$}}}) \r)\circ C[\eth^{\mbox{{\tiny$W$}}}] \oplus C[\eth^{\mbox{{\tiny$W\pr$}}}].\nonumber  
\end{eqnarray}

Next we show that the diagram has exact rows and is commutative up to adding a smoothing operator to the vertical Fredholm maps. We may write 
relative to \eqref{idents} and using \lemref{C-Pi orientation} 
$$
\Pi_>^{\mbox{{\tiny $\pd W$}}}\oplus U^{\mbox{{\tiny $\pd W$}}}
 = \left(\begin{array}{cc}
\Pi_<^{\mbox{{\tiny $X$}}}\oplus U_-^{\mbox{{\tiny $X$}}}  & 0 \\
 0 & \Pi_>^{\mbox{{\tiny $N$}}}\oplus U_+^{\mbox{{\tiny $N$}}}
 \end{array}\right)  \hf \in \  \Psi^0(X\sqcup N)$$
 with $U_+^{\mbox{{\tiny $X$}}} = \Pi^{ev}_0(B_X)$ and  $U_-^{\mbox{{\tiny $X$}}} = \Pi^{odd}_0(B_X)$,  mindful of \lemref{B even odd}.  While
 $$
\Pi_>^{\mbox{{\tiny $\pd W\pr$}}}\oplus U^{\mbox{{\tiny $\pd W\pr$}}}
 = \left(\begin{array}{cc}
\Pi_<^{\mbox{{\tiny $N$}}}\oplus U_-^{\mbox{{\tiny $N$}}}  & 0 \\
 0 & \Pi_>^{\mbox{{\tiny $Y$}}}\oplus U_+^{\mbox{{\tiny $Y$}}}
 \end{array}\right)  \hf \in \  \Psi^0(N\sqcup Y),$$
 $$
\Pi_>^{\mbox{{\tiny $\pd(W \cup_\phi W\pr)$}}}\oplus U^{\mbox{{\tiny $\pd(W \cup_\phi W\pr)$}}}
 = \left(\begin{array}{cc}
\Pi_<^{\mbox{{\tiny $X$}}}\oplus U_-^{\mbox{{\tiny $X$}}}  & 0 \\
 0 & \Pi_>^{\mbox{{\tiny $Y$}}}\oplus U_+^{\mbox{{\tiny $Y$}}}
 \end{array}\right)  \hf \in \  \Psi^0(X\sqcup Y).$$
These choices for the projections $U^V_\pm$ provide a canonical identification
$$ \ran(\Pi_>^{\mbox{{\tiny $\pd(W \cup_\phi W\pr) $}}}\oplus U^{\mbox{{\tiny $\pd(W \cup_\phi W\pr) $}}})  =  \ran(\Pi_<^{\mbox{{\tiny $X$}}}\oplus U_-^{\mbox{{\tiny $X$}}}) \oplus   \ran(\Pi_>^{\mbox{{\tiny $Y$}}}\oplus U_+^{\mbox{{\tiny $Y$}}})$$
and, since 
$ (\Pi_>^{\mbox{{\tiny $N$}}}\oplus U_+^{\mbox{{\tiny $N$}}}) \oplus (\Pi_<^{\mbox{{\tiny $N$}}}\oplus U_-^{\mbox{{\tiny $N$}}})  = id_N,$
a canonical identification
\begin{equation}\label{Mayer Vietoris 5}  
 \ran(\Pi_>^{\mbox{{\tiny $\pd W$}}}\oplus U^{\mbox{{\tiny $\pd W$}}}) \oplus \ran(\Pi_>^{\mbox{{\tiny $\pd W\pr$}}}\oplus U^{\mbox{{\tiny $\pd W\pr$}}})  =  \ran(\Pi_<^{\mbox{{\tiny $X$}}}\oplus U_-^{\mbox{{\tiny $X$}}}) \oplus H_N \oplus  \ran(\Pi_>^{\mbox{{\tiny $Y$}}}\oplus U_+^{\mbox{{\tiny $Y$}}}),
\end{equation}
hence defining the maps in the lower exact sequence of the diagram.  

The exactness of the top row has been accounted for above. As  $K(\eth^{\mbox{{\tiny $W \cup_\phi W\pr $}}}) \< H_X \oplus H_Y$, an element $\z\in K(\eth^{\mbox{{\tiny $W \cup_\phi W\pr $}}})$ may be written uniquely as  
$\z = (\xi_X, \eta_Y)$ with $ \xi_X\in H_X, \ \eta_Y\in H_Y.$
For convenience, and since it does not affect any previous construction, we also include the involution $(\a,\b) \mto (\a,-\b)$ on $K(\eth^{\mbox{{\tiny $W\pr $}}}) \< H_N \oplus H_Y$, so that the inclusion $$K(\eth^{\mbox{{\tiny $W \cup_\phi W\pr $}}}) \to K(\eth^{\mbox{{\tiny $W$}}})\oplus K(\eth^{\mbox{{\tiny $W\pr$}}}) \hf  \iis \ \ 
(\xi_X, \eta_Y) \mto (\xi_X,\nu_N)\oplus (-\nu_N,\eta_Y),$$
where
$\nu_N = \nu_N(\xi_X, \eta_Y)$
is uniquely defined via unique continuation and the Poisson operator; $(\xi_X, \eta_Y)$ corresponds uniquely via the Poisson operator  to an element of $\Ker(\eth^{\mbox{{\tiny $W \cup_\phi W\pr $}}})$, then restrict to the hypersurfaces $X$, $N$ and $Y$.  

Now replace $G_1$ by 
$
\Gg_1 =((\Pi_<^{\mbox{{\tiny $X$}}}\oplus U_-^{\mbox{{\tiny $X$}}})\oplus I_N)\circ C[\eth^{\mbox{{\tiny$W$}}}] + 
(I_N\oplus (\Pi_>^{\mbox{{\tiny $Y$}}}\oplus U_+^{\mbox{{\tiny $Y$}}}))\circ C[\eth^{\mbox{{\tiny$W\pr$}}}]$
as a map 
$$K(\eth^{\mbox{{\tiny $W$}}})\oplus K(\eth^{\mbox{{\tiny $W\pr$}}}) \to \ran(\Pi_<^{\mbox{{\tiny $X$}}}\oplus U_-^{\mbox{{\tiny $X$}}}) \oplus H_N \oplus  \ran(\Pi_>^{\mbox{{\tiny $Y$}}}\oplus U_+^{\mbox{{\tiny $Y$}}}),$$
where $C[\eth^{\mbox{{\tiny$W$}}}]$ and $(\Pi_<^{\mbox{{\tiny $X$}}}\oplus U_-^{\mbox{{\tiny $X$}}})\oplus I_N$  mean $C[\eth^{\mbox{{\tiny$W$}}}]\oplus 0$ and $(\Pi_<^{\mbox{{\tiny $X$}}}\oplus U_-^{\mbox{{\tiny $X$}}})\oplus I_N\oplus 0$, and so on. 

\begin{lem}\label{diagram commutes}
With $G_1$ replaced by $\Gg_1$ the diagram \eqref{diagram} commutes.
\end{lem}
Proof:\;  
 $\Gg_1$ evaluated on  $(\xi_X, \la_N) \oplus (\mu_N, \eta_Y) \in K(\eth^{\mbox{{\tiny $W$}}})\oplus K(\eth^{\mbox{{\tiny $W\pr$}}})$ is 
$
\Gg_1((\xi_X, \la_N) , (\mu_N, \eta_Y)) = ((\Pi_<^{\mbox{{\tiny $X$}}}\oplus U_-^{\mbox{{\tiny $X$}}})\xi_X, \ \la_N + \mu_N, \ (\Pi_>^{\mbox{{\tiny $Y$}}}\oplus U_+^{\mbox{{\tiny $Y$}}})\eta_Y).$
With $G_1$ replaced by $\Gg_1$: the left-hand square of  \eqref{diagram} is 
$$
\begin{array}{ccc}
(\xi_X, \eta_Y)  & \to  &  ((\xi_X,\la), (-\la, \eta_Y)) \\[2mm]
 \downarrow  &   &   \downarrow   \\[3mm]
 \hf ((\Pi_<^{\mbox{{\tiny $X$}}}\oplus U_-^{\mbox{{\tiny $X$}}})\xi_X, (\Pi_>^{\mbox{{\tiny $Y$}}}\oplus U_+^{\mbox{{\tiny $Y$}}}) \eta_Y)  & \to  &  \hf ((\Pi_<^{\mbox{{\tiny $X$}}}\oplus U_-^{\mbox{{\tiny $X$}}})\xi_X, \, 0, \,  (\Pi_>^{\mbox{{\tiny $Y$}}}\oplus U_+^{\mbox{{\tiny $Y$}}}) \eta_Y)  
 \end{array}
$$
and the right-hand square is 
$$
\begin{array}{ccc}
 ((\xi_X,\la_N ), (\mu_N, \eta_Y))  & \to  & \la_N  + \mu_N \\[2mm]
 \downarrow  &   &   \downarrow   \\[3mm]
 \hf ((\Pi_<^{\mbox{{\tiny $X$}}}\oplus U_-^{\mbox{{\tiny $X$}}})\xi_X, \, \la_N  + \mu_N, \, (\Pi_>^{\mbox{{\tiny $Y$}}}\oplus U_+^{\mbox{{\tiny $Y$}}}) \eta_Y)  & \to  &  \la_N  + \mu_N.
 \end{array}
$$
\begin{flushright}
$\Box$
\end{flushright} 

\begin{lem}\label{G - Gg}  
$$G_1 - \Gg_1 : K(\eth^{\mbox{{\tiny $W$}}})\oplus K(\eth^{\mbox{{\tiny $W\pr$}}}) \to \ran(\Pi_<^{\mbox{{\tiny $X$}}}\oplus U_-^{\mbox{{\tiny $X$}}}) \oplus H_N \oplus  \ran(\Pi_>^{\mbox{{\tiny $Y$}}}\oplus U_+^{\mbox{{\tiny $Y$}}})$$
is the restriction of a smoothing operator $H_X \oplus H_N \oplus H_N \oplus H_Y\to H_X \oplus H_N\oplus H_Y$.
\end{lem}
Proof:\;  For $(\xi_X, \la_N) \oplus (\mu_N, \eta_Y) \in K(\eth^{\mbox{{\tiny $W$}}})\oplus K(\eth^{\mbox{{\tiny $W\pr$}}})$ 
$$G_1((\xi_X, \la_N) , (\mu_N, \eta_Y)) := ((\Pi_<^{\mbox{{\tiny $X$}}}\oplus U_-^{\mbox{{\tiny $X$}}})\xi_X, \ (\Pi_>^{\mbox{{\tiny $N$}}}\oplus U_+^{\mbox{{\tiny $N$}}})\la_N + (\Pi_<^{\mbox{{\tiny $N$}}}\oplus U_-^{\mbox{{\tiny $N$}}})\mu_N, \ (\Pi_>^{\mbox{{\tiny $Y$}}}\oplus U_+^{\mbox{{\tiny $Y$}}})\eta_Y).$$
Hence
$(G_1 - \Gg_1)((\xi_X, \la_N) , (\mu_N, \eta_Y)) = (0,  \, (\Pi_<^{\mbox{{\tiny $N$}}}\oplus U_-^{\mbox{{\tiny $N$}}})\la_N + (\Pi_<^{\mbox{{\tiny $N$}}}\oplus U_+^{\mbox{{\tiny $N$}}})\mu_N, \, 0).$
Since $U_\pm^{\mbox{{\tiny $N$}}}$ is smoothing we may ignore this term, and it is enough to show that
$
(\xi_X, \la_N) \to  (0,  \, \Pi_<^{\mbox{{\tiny $N$}}}\la_N)\hf \aand\hf (\mu_N, \eta_Y) \to (\Pi_<^{\mbox{{\tiny $N$}}}\mu_N, \, 0)
$
are (restrictions of) smoothing operators. For this, on $(\xi_X, \la_N) \in K(\eth^{\mbox{{\tiny $W$}}})= \ran(C[\eth^{\mbox{{\tiny$W$}}}](\xi_X, \la_N) )$ one has 
$(\xi_X, \la_N) = C[\eth^{\mbox{{\tiny$W$}}}](\xi_X, \la_N)$. Writing $C[\eth^{\mbox{{\tiny$W$}}}] = \left(\begin{array}{cc} C^{\mbox{{\tiny $X,X$}}} & C^{\mbox{{\tiny $N,X$}}}  \\ C^{\mbox{{\tiny $X,N$}}}  & C^{\mbox{{\tiny $N,N$}}} \end{array}\right)$ as a 2x2 block matrix on $H_X\oplus H_N$, we see $C^{\mbox{{\tiny $X,N$}}} : H_X \to H_N$ and $C^{\mbox{{\tiny $N,X$}}} : H_N\to H_X$ are smoothing, in view of \eqref{APS - C}, this gives
$\la_N =  C^{\mbox{{\tiny $X,N$}}}\xi_X  + C^{\mbox{{\tiny $N,N$}}}\la_N$
and that the first of the maps in question is the restriction of 
$\left(\begin{array}{cc} 0 & 0\\ 
\Pi_<^{\mbox{{\tiny $N$}}}C^{\mbox{{\tiny $X,N$}}}  & \Pi_<^{\mbox{{\tiny $N$}}}C^{\mbox{{\tiny $N,N$}}} \end{array}\right) \in \Psi^\Z(X\sqcup N).$
Since $C^{\mbox{{\tiny $X,N$}}}$ is smoothing, we have only to show that $\Pi_<^{\mbox{{\tiny $N$}}}C^{\mbox{{\tiny $N,N$}}}\in \Psi^{-\oo}(N).$ But \eqref{APS - C} states
$\left(\begin{array}{cc} C^{\mbox{{\tiny $X,X$}}} & C^{\mbox{{\tiny $N,X$}}}  \\ C^{\mbox{{\tiny $X,N$}}}  & C^{\mbox{{\tiny $N,N$}}} \end{array}\right)  -  \left(\begin{array}{cc} \Pi_<^{\mbox{{\tiny $X$}}} & 0  \\ 0 & \Pi_>^{\mbox{{\tiny $N$}}}\end{array}\right)   \ \in \Psi^{-\oo}(X\sqcup N)$
and, in particular, that $C^{\mbox{{\tiny $N,N$}}} - \Pi_>^{\mbox{{\tiny $N$}}} \in \Psi^{-\oo}(N).$ Hence, 
$\Pi_<^{\mbox{{\tiny $N$}}}C^{\mbox{{\tiny $N,N$}}} = \Pi_<^{\mbox{{\tiny $N$}}}(C^{\mbox{{\tiny $N,N$}}} - \Pi_>^{\mbox{{\tiny $N$}}} )$
is smoothing. 
\hfill
$\Box$\\[2mm]

Since $G_1$ is from \eqref{G1a} the direct sum of the operators $(\Pi_>^{\mbox{{\tiny $\pd W$}}}\oplus U^{\mbox{{\tiny $\pd W$}}})\circ C[\eth^{\mbox{{\tiny$W$}}}] : K(\eth^{\mbox{{\tiny $W$}}}) \to \ran(\Pi_>^{\mbox{{\tiny $\pd W $}}}\oplus U^{\mbox{{\tiny $\pd W $}}}) $ and 
$(\Pi_>^{\mbox{{\tiny $\pd W\pr$}}}\oplus U^{\mbox{{\tiny $\pd W\pr$}}})\circ C[\eth^{\mbox{{\tiny$W\pr$}}}] : K(\eth^{\mbox{{\tiny $W\pr$}}}) \to \ran(\Pi_>^{\mbox{{\tiny $\pd W\pr $}}}\oplus U^{\mbox{{\tiny $\pd W\pr $}}}) $
and from \eqref{ind PK} these are Fredholm, then $G_1$ is a Fredholm operator with index $$\ind(G_1) = \ind \left( (\Pi_>^{\mbox{{\tiny $\pd W$}}}\oplus U^{\mbox{{\tiny $\pd W$}}})\circ C[\eth^{\mbox{{\tiny$W$}}}]\right)  \ + \ \ind \left((\Pi_>^{\mbox{{\tiny $\pd W\pr$}}}\oplus U^{\mbox{{\tiny $\pd W\pr$}}})\circ C[\eth^{\mbox{{\tiny$W\pr$}}}]\right).$$ 
By \lemref{G - Gg}  
$\ind(G_1) \ = \  \ind(\Gg_1).$
By \lemref{diagram commutes} and Lemma 5 on p.202 of \cite{Mac} 
$ \ind(\Gg_1) = \ind(G_0)  + \ind(id_{H_{\mbox{{\tiny $N$}}}})= \ind(G_0).$
Hence
$\ind(G_0) \ = \  \ind(G_1).$
That is, 
$\Pi_>^{\mbox{{\tiny $\pd(W \cup_\phi W\pr) $}}}\oplus U^{\mbox{{\tiny $\pd(W \cup_\phi W\pr)$}}})\circ C[\eth^{\mbox{{\tiny$W_0 \cup_\phi W_1$}}}]$ has index equal to $   \ind \left( (\Pi_>^{\mbox{{\tiny $\pd W$}}}\oplus U^{\mbox{{\tiny $\pd W$}}})\circ C[\eth^{\mbox{{\tiny$W$}}}]\right)  \ + \ \ind \left((\Pi_>^{\mbox{{\tiny $\pd W\pr$}}}\oplus U^{\mbox{{\tiny $\pd W\pr$}}})\circ C[\eth^{\mbox{{\tiny$W\pr$}}}]\right).$
But
\begin{eqnarray*}
\ind \left( (\Pi_>^{\mbox{{\tiny $\pd W$}}}\oplus U^{\mbox{{\tiny $\pd W$}}})\circ C[\eth^{\mbox{{\tiny$W$}}}]\right)  & \stackrel{\eqref{ind=tr} }{=} & \Tr_{\mbox{{\tiny $X\sqcup N$}}}\l(C[\eth^{\mbox{{\tiny$W$}}}] -   \Pi_>^{\mbox{{\tiny $\pd W$}}}\oplus U^{\mbox{{\tiny $\pd W$}}}\r) \\[2mm] 
& \stackrel{\eqref{commuting trace reduced}}{=} & \wt\Tr_{\mbox{{\tiny $X\sqcup N$}}}\l( \pi_{\mbox{{\tiny $X\sqcup N$}}} \l(C[\eth^{\mbox{{\tiny$W$}}}] -   \Pi_>^{\mbox{{\tiny $\pd W$}}}\oplus U^{\mbox{{\tiny $\pd W$}}}\r)\r) \\[2mm] 
&  \stackrel{\eqref{tilde tracial Finfty}}{=} & \wt\Tr_{\mbox{{\tiny $X\sqcup N \sqcup N \sqcup Y$}}}\l( \wt\eta_{\mbox{{\tiny $N\sqcup Y$}}} \l(\pi_{\mbox{{\tiny $X\sqcup N$}}} \l( C[\eth^{\mbox{{\tiny$W$}}}] -   \Pi_>^{\mbox{{\tiny $\pd W$}}}\oplus U^{\mbox{{\tiny $\pd W$}}}\r)\r)\r),
\end{eqnarray*}
\begin{eqnarray*}
\ind \left((\Pi_>^{\mbox{{\tiny $\pd W\pr$}}}\oplus U^{\mbox{{\tiny $\pd W\pr$}}})\circ C[\eth^{\mbox{{\tiny$W\pr$}}}]\right) & = &   \wt\Tr_{\mbox{{\tiny $X\sqcup N \sqcup N \sqcup Y$}}}\l( \wt\eta_{\mbox{{\tiny $X\sqcup N$}}} \l(\pi_{\mbox{{\tiny $X\sqcup N$}}} \l( C[\eth^{\mbox{{\tiny$W\pr$}}}] -   \Pi_>^{\mbox{{\tiny $\pd W\pr$}}}\oplus U^{\mbox{{\tiny $\pd W\pr$}}}\r)\r)\r),
\end{eqnarray*}
$$\ind \left( (\Pi_>^{\mbox{{\tiny $\pd(W \cup_\phi W\pr) $}}}\oplus U^{\mbox{{\tiny $\pd(W \cup_\phi W\pr)$}}})\circ C[\eth^{\mbox{{\tiny$W_0 \cup_\phi W_1$}}}]\right) $$
$$ =  \wt\Tr_{\mbox{{\tiny $X\sqcup N \sqcup N \sqcup Y$}}}\l( \wt\eta_{\mbox{{\tiny $N\sqcup N$}}} \l(\pi_{\mbox{{\tiny $X\sqcup Y$}}} \l( C[\eth^{\mbox{{\tiny$W_0 \cup_\phi W_1$}}}] -   \Pi_>^{\mbox{{\tiny $\pd(W \cup_\phi W\pr) $}}}\oplus U^{\mbox{{\tiny $\pd(W \cup_\phi W\pr)$}}}\r)\r)\r).$$
The (reduced) trace $\wt\Tr_{\mbox{{\tiny $X\sqcup N \sqcup N \sqcup Y$}}}$ therefore vanishes on the element
$$\wt\eta_{\mbox{{\tiny $N\sqcup N$}}} \l(\pi_{\mbox{{\tiny $X\sqcup Y$}}} \l( C[\eth^{\mbox{{\tiny$W_0 \cup_\phi W_1$}}}] -   \Pi_>^{\mbox{{\tiny $\pd(W \cup_\phi W\pr) $}}}\oplus U^{\mbox{{\tiny $\pd(W \cup_\phi W\pr)$}}}\r)\r) $$ $$ -  \wt\eta_{\mbox{{\tiny $X\sqcup N$}}} \l(\pi_{\mbox{{\tiny $X\sqcup N$}}} \l( C[\eth^{\mbox{{\tiny$W\pr$}}}] -   \Pi_>^{\mbox{{\tiny $\pd W\pr$}}}\oplus U^{\mbox{{\tiny $\pd W\pr$}}}\r)\r) -\wt\eta_{\mbox{{\tiny $N\sqcup Y$}}} \l(\pi_{\mbox{{\tiny $X\sqcup N$}}} \l( C[\eth^{\mbox{{\tiny$W$}}}] -   \Pi_>^{\mbox{{\tiny $\pd W$}}}\oplus U^{\mbox{{\tiny $\pd W$}}}\r)\r)$$
in
$\frac{\Fs_{-\oo}(X\sqcup N \sqcup N\sqcup Y)}{[\Fs_{-\oo}(X\sqcup N \sqcup N\sqcup Y), \Fs_{-\oo}(X\sqcup N \sqcup N\sqcup Y)]}$
By \eqref{trace reduced}, this element is zero, which is \eqref{log sgn add 3}. 
 $\hfill\Box$

\vskip 5mm

A closer look at the identity \eqref{log sgn add} reveals that it is equivalent to the Calderon projections fitting together with respect to gluing in the following way:
\begin{cor}\label{additivity signature log proj}  
With $C(\eth^{\mbox{{\tiny$W_1$}}})^\bot : = (I \oplus 0) - C(\eth^{\mbox{{\tiny$W_1$}}}) \in \Psi^0(M_1\sqcup M_2).$, one has 
$$\eta_{\mbox{{\tiny $M_1$}}}C(\eth^{\mbox{{\tiny$W_0 \cup_{M_1}W_1$}}})  - \eta_{\mbox{{\tiny $M_2$}}}C(\eth^{\mbox{{\tiny$W_0$}}}) - \eta_{\mbox{{\tiny $M_0$}}}C(\eth^{\mbox{{\tiny$W_1$}}})^\bot \in [\Fs_{\mbox{{\tiny $-\oo$}}}(M_0\sqcup M_1\sqcup M_2), \Fs_{\mbox{{\tiny $-\oo$}}}(M_0\sqcup M_1\sqcup M_2)].$$
\end{cor}

\vskip 10mm

{\small \noi \textsc{Department of Mathematics \\King's College London}

\end{document}